\documentclass[12pt,reqno,a4paper]{amsart}

\usepackage{amsmath,amsfonts,amsthm,amssymb}
\usepackage{cite}

\usepackage{a4wide}

\usepackage[dvips]{hyperref}

\usepackage[footnotesize,hang,sf]{caption}

\usepackage{dcolumn}

\usepackage{txfonts,times}
\newcommand{\widebar}{\overline}
\usepackage{titlesec}
\titleformat{\section}
{\bfseries\scshape\centering}
{\thesection.}{.5em}{}
\titleformat{\subsection}
{\rmfamily\bfseries}
{\thesubsection}{.5em}{}
\usepackage{,nicefrac}
\newcommand{\vv}{\vec}
\usepackage{psfrag}
\usepackage[dvips]{graphicx}
\numberwithin{equation}{section}

\newtheorem{theorem}{Theorem}
\newtheorem{prop}[theorem]{Proposition}
\newtheorem{coro}[theorem]{Corollary}
\newtheorem{lemma}[theorem]{Lemma}
\newtheorem{conj}{Conjecture}
\DeclareMathOperator{\sign}{sign}
\DeclareMathOperator{\Tr}{Tr}

\DeclareMathOperator{\CS}{CS}

\DeclareMathOperator{\Vol}{Vol}
\DeclareMathOperator*{\Res}{Res}

\DeclareMathOperator{\eveodd}{m_2}
\newcommand{\I}{\mathrm{i}}
\newcommand{\E}{\mathrm{e}}

\begin{document}

\renewcommand{\thefootnote}{\fnsymbol{footnote}}
\baselineskip 16pt
\parskip 6.4pt
\sloppy




\title{Quantum Invariants, Modular Forms, and Lattice Points II
}


    \author{Kazuhiro \textsc{Hikami}}


  \address{Department of Physics, Graduate School of Science,
    University of Tokyo,
    Hongo 7--3--1, Bunkyo, Tokyo 113--0033,   Japan.
    }

%
    \urladdr{\url{http://gogh.phys.s.u-tokyo.ac.jp/~hikami/}}

    \email{\texttt{hikami@phys.s.u-tokyo.ac.jp}}


\date{February 24, 2006}

\begin{abstract}
We study the SU(2) Witten--Reshetikhin--Turaev invariant for the
Seifert fibered homology spheres with $M$-exceptional fibers.
We show that the WRT invariant can be written in terms of 
(differential of) the Eichler
integrals of modular forms with weight $1/2$ and
$3/2$.
By use of nearly modular property of the Eichler integrals we shall
obtain asymptotic expansions of the WRT invariant in the large-$N$
limit.
We further reveal that the number of the gauge equivalent classes of
flat
connections, which dominate the asymptotics  of the WRT
invariant in $N\to\infty$,
is related to the number of integral
lattice points inside the $M$-dimensional tetrahedron.

\end{abstract}





\maketitle

\section{Introduction}

The Witten  invariant for the
3-manifold $\mathcal{M}$ is defined by the Chern--Simons path integral
as~\cite{EWitt89a}
(see also Ref.~\citen{Atiya90Book})
\begin{equation}
  \label{define_Witten_partition}
  Z_k(\mathcal{M})
  =
  \int \exp \left(
    2 \, \pi \, \I \, k \, \CS(A)
  \right) \,
  \mathcal{D}A
\end{equation}
where $k \in \mathbb{Z}$, and $\CS(A)$ is the Chern--Simons functional
\begin{equation}
  \label{CS_functional}
  \CS(A)
  =
  \frac{1}{8 \, \pi^2} \int_\mathcal{M}
  \Tr
  \left(
    A \wedge \mathrm{d} A
    + 
    \frac{2}{3} \, A \wedge A \wedge A
  \right)
\end{equation}
In a limit $k\to\infty$ of  the Witten invariant $Z_k(\mathcal{M})$, we
may apply the saddle point method.
As the saddle point of the Chern--Simons
functional~\eqref{CS_functional} denotes the flat connection
\begin{equation}
\mathrm{d} A + A \wedge A=0
\end{equation}
the asymptotics of the partition function becomes
a sum of the Chern--Simons invariants,
and it is expected to be~~\cite{EWitt89a,Atiya90Book,FreeGomp91a}
\begin{multline}
  \label{expected_asymptotics}
  Z_k(\mathcal{M}) 
  \sim
  \frac{1}{2} \, \E^{-\frac{3}{4} \pi \I ( 1 + b^1)}
  \sum_\alpha
  (k+2)^{(\dim H^1 - \dim H^0)/2}
  \\
  \times
  \sqrt{T_\alpha} \,
  \E^{- 2 \pi \I (I_\alpha/4 + \dim H^0 /8)}
  \,
  \E^{2 \pi \I (k+2) \CS(A_\alpha)} \,
\end{multline}
Here the sum of $\alpha$ denotes a gauge equivalent class of flat
connections,
and $T_\alpha$ and
$I_\alpha$ respectively denote
the Reidemeister torsion and the spectral
flow.
The first Betti number is $b^1$, and $H^i$ is the cohomology space.

To study the asymptotic behavior of the Witten invariant 
rigorously, we need explicit expression of the invariant.
Alternative and  combinatorial
definition  of this quantum invariant
was given by Reshetikhin and
Turaev~\cite{ResheTurae91a}
(see also Ref.~\citen{KirbMelv91a}).
We denote
$\tau_N(\mathcal{M})$
as  the Witten--Reshetikhin--Turaev (WRT) invariant, which is related to
the Witten invariant $Z_k(\mathcal{M})$ by
\begin{equation}
  Z_k(\mathcal{M})
  =
  \frac{
    \tau_{k+2} \left( \mathcal{M} \right)
  }{
    \tau_{k+2} \left( S^2 \times S^1 \right)
  }
\end{equation}
and we have
\begin{gather*}
  \tau_N \left( S^3 \right) =1
  \\[2mm]
  \tau_N \left( S^2 \times S^1 \right)
  =
  \sqrt{\frac{N}{2}} \,
  \frac{1}{
    \sin(\pi/N)
  }
\end{gather*}
Using this definition of the WRT invariant, asymptotic behavior of the
WRT invariants for certain 3-manifolds has been extensively
studied~\cite{LawreRozan99a,RLawre95a,Rozan94,Rozan94b,Rozan95a,Rozan96c,Rozan96d,SKHanse05a}.

Several years ago,
Lawrence and Zagier found a connection between the WRT invariant and
 modular form~\cite{LawrZagi99a}.
They showed  that the  WRT invariant $\tau_N(\mathcal{M})$ for the Poincar{\'e} homology
sphere $\mathcal{M}=\Sigma(2,3,5)$ can be regarded as a limiting value of the Eichler integral of
vector modular form with weight $3/2$.
Thanks to this correspondence, 
the exact asymptotic
expansion of the WRT invariant in the large-$N$ limit can be computed,
and
topological invariants such as the Chern--Simons invariant and the
Reidemeister torsion can be interpreted from the viewpoint of  modular forms.
Meanwhile it has   been established
that this remarkable structure of the quantum invariants
holds for
the WRT invariants for 3-manifolds such as
the Brieskorn homology spheres~\cite{KHikami04b},
4-exceptional fibered Seifert homology spheres~\cite{KHikami04e},
and the spherical Seifert
manifolds~\cite{KHikami05a}.
Also established is a connection between the Eichler integrals of
vector modular forms with weight $1/2$ and
the special values of the colored Jones polynomial for the torus
knot $\mathcal{T}_{s,t}$~\cite{KHikami03c}
(see also Refs.~\citen{KHikami02b,KHikami03b,KHikami02c,DZagie01a})
and the torus link $\mathcal{T}_{2,2m}$~\cite{KHikami03a}.

One of the benefit of the quantum invariant/modular form
correspondence is  an observation
that a limiting value of the Ramanujan mock theta functions~\cite{Ramanujan87Book}
in $q\to\E^{2 \pi \I /N}$ from outside a unit circle
coincides with the
WRT invariants for the spherical Seifert manifolds~\cite{KHikami05b}.
This fact opens up a new insight to  modular forms
and the Ramanujan mock theta functions,
and we can expect that further
studies on the quantum invariant/modular form correspondence 
should be fruitful.

In this article, as a continuation of
Refs.~\citen{KHikami04b,KHikami04e},  we study an exact asymptotic
expansion of the
WRT invariant
$\tau_N(\mathcal{M})$
for the $M$-exceptional fibered Seifert
integral homology sphere
$\mathcal{M}=\Sigma ( p_1, p_2, \dots, p_M)$,
where $p_j$ are pairwise
coprime positive integers.
By use of  modular forms with half-integral weight, we derive an
asymptotic expansion in $N\to\infty$ number theoretically.

This paper is organized as follows.
In Section~\ref{sec:WRT} we review the construction of the WRT
invariant for the Seifert fibered homology spheres following
Ref.~\citen{LawreRozan99a}.
An explicit form of the WRT invariant is given.
Also discussed is an integral expression of the invariant.
In Section~\ref{sec:modular_form} we introduce a family of vector
modular forms with half-integral weight.
We define the Eichler integrals thereof, and study the nearly modular
property of a limiting value of the Eichler integrals.
By use of this quasi modular transformation property,
we compute the asymptotic expansion of the WRT invariant in the
large-$N$ limit in Section~\ref{sec:asymptotic}.
We shall see that the invariant is a limiting value of the holomorphic
function~\cite{KHabi02a}.
We study  a contribution of dominating terms in the large-$N$ limit in
detail,
and reveal a relationship with the number of the integral lattice
points inside the higher dimensional tetrahedron.
Also given is
an explicit relationship between the Casson invariant and 
the first non-trivial coefficient of
the Ehrhart polynomial.
In  Section~\ref{sec:examples} we  give some results based on
numerical computations.
We compare the exact value of the WRT invariant with our asymptotic formula.
The last section is devoted to conclusion and discussions.

\section{WRT Invariant for Seifert Integral Homology Sphere}
\label{sec:WRT}

\begin{figure}
  \centering
  
  \begin{psfrags}
    \psfrag{x}{$0$}
    \psfrag{a}{$\nicefrac{p_1}{q_1}$}
    \psfrag{b}{$\nicefrac{p_2}{q_2}$}
    \psfrag{c}{$\nicefrac{p_3}{q_3}$}
    \psfrag{d}{$\nicefrac{p_M}{q_M}$}
    \includegraphics[scale=0.36, bb=-240 -280 80 320]{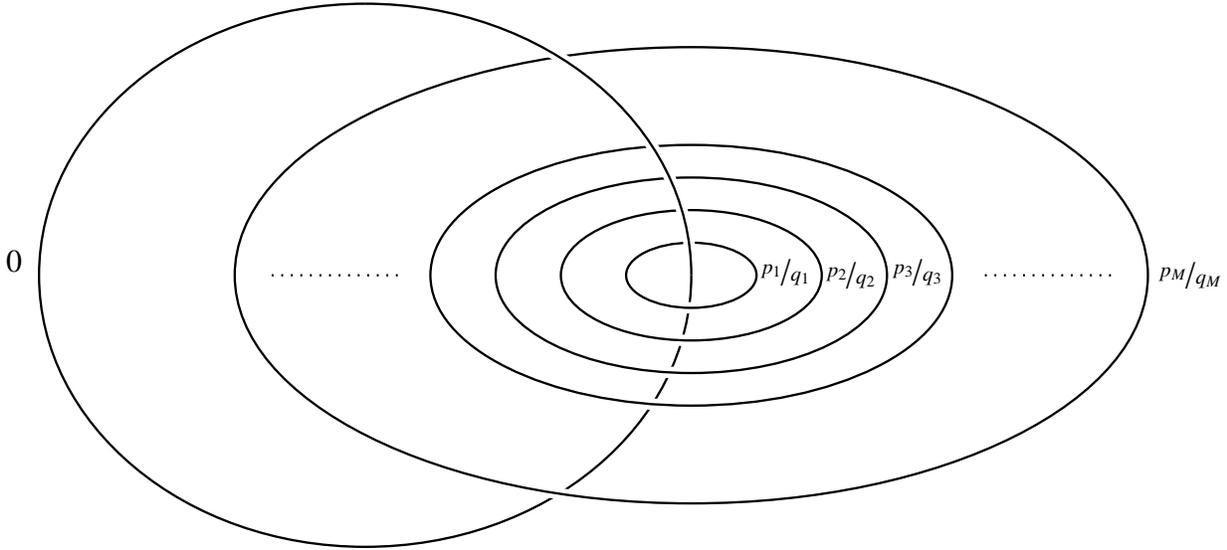}    
  \end{psfrags}
  \caption{Surgery description of the Seifert homology sphere $\Sigma(p_1,\dots,p_M)$}
  \label{fig:surgery}
\end{figure}

Following Ref.~\citen{LawreRozan99a},
we compute the WRT
invariant $\tau_N(\mathcal{M})$ for the
Seifert fibered integral  homology sphere with $M$-exceptional fibers
$
\mathcal{M}=\Sigma(\vv{p})
=\Sigma(p_1, p_2, \dots, p_M)
$
where $p_j$ are pairwise coprime positive integers.
Hereafter we use
$\vv{p}$ as $M$-tuple
\begin{equation*}
  \vv{p} =
  (p_1, p_2, \dots, p_M)
\end{equation*}
The Seifert fibered integral homology sphere $\Sigma(\vv{p})$
has a rational surgery description as
Fig.~\ref{fig:surgery}
(see, \emph{e.g.}, Refs.~\citen{SeifeThrel80,Orlik72,NSavel02Book}),
and the fundamental group has a   presentation
\begin{equation}
  \label{fundamental_group}
  \pi_1
  \left(
    \Sigma(\vv{p})
  \right)
  =
  \left\langle
    x_1, x_2, \dots, x_M , h
    ~
    \Big|
    ~
    \begin{array}{c}
      \text{$h$ is center}
      \\[1mm]
      \text{$x_j^{~p_j} = h^{-q_j}$
        for $1 \leq j \leq M$}
      \\[1mm]
      x_1 \, x_2 \cdots x_M=1
    \end{array}
  \right\rangle
\end{equation}
Here $q_j \in \mathbb{Z}$ is coprime to $p_j$, and we
have a constraint so that the fundamental
group~\eqref{fundamental_group} gives the homology sphere;
\begin{equation}
  P \sum_{j=1}^M \frac{q_j}{p_j} =1
\end{equation}
Here and hereafter we use
\begin{equation}
  P =  P(\vv{p})=
  \prod_{j=1}^M p_j
\end{equation}

When the 3-manifold $\mathcal{M}$ is constructed by the rational
surgeries
$p_j / q_j$
on the $j$-th component of $n$-component link $\mathcal{L}$, it was
shown~\cite{ResheTurae91a,LCJeff92a} that the SU(2) WRT
invariant $\tau_N(\mathcal{M})$ is
given by
\begin{equation}
  \label{define_invariant}
  \tau_N(\mathcal{M})
  =
  \mathrm{e}^{
    \frac{\pi \mathrm{i}}{4} \frac{N-2}{N}
    \left(
      \sum_{j=1}^n \Phi(U^{(p_j,q_j)})
      -3 \sign(\mathbf{L})
    \right)
  }
  \sum_{k_1, \dots, k_n=1}^{N-1}
  J_{k_1,\dots,  k_n}(\mathcal{L}) \,
  \prod_{j=1}^n
  \rho(U^{(p_j,q_j)})_{k_j, 1}
\end{equation}
Here the surgery data $p_j/q_j$  is encoded  by an $SL(2; \mathbb{Z})$ matrix
\begin{equation*}
  U^{(p_j,q_j)}=
  \begin{pmatrix}
    p_j & r_j \\
    q_j & s_j
  \end{pmatrix}
\end{equation*}
The Rademacher $\Phi$-function $\Phi(U)$ is defined by~\cite{HRadema73}
\begin{equation}
  \Phi
  \left(
    \begin{pmatrix}
      p & r \\
      q & s
    \end{pmatrix}
  \right)
  =
  \begin{cases}
    \displaystyle
    \frac{p+s}{q} - 12 \, s(p, q) &
    \text{for $q \neq 0$}
    \\[4.8mm]
    \displaystyle
    \frac{r}{s}
    &
    \text{for $q=0$}
  \end{cases}
\end{equation}
where   $s(b,a)$ denotes the Dedekind sum~\eqref{Dedekind_sum}.
An $n\times n$ matrix $\mathbf{L}$ is a linking matrix
\begin{equation}
  \mathbf{L}_{j,k} = \mathrm{lk}(j,k) + \frac{p_j}{q_j} \cdot \delta_{j,k}  
\end{equation}
where $\mathrm{lk}(j,k)$ denotes the linking number of the $j$- and
$k$-th components of link $\mathcal{L}$,
and $\sign(\mathbf{L})$ denotes a signature of $\mathbf{L}$, \emph{i.e.},
the
difference between the number of positive and negative eigenvalues of
$\mathbf{L}$.
The  polynomial $J_{k_1,\dots,k_n}(\mathcal{L})$ denotes the colored Jone
polynomial for link $\mathcal{L}$ with  color $k_j$ for the $j$-th
component link, and
$\rho(U^{(p,q)})$ 
is a representation $\rho$ of $PSL(2;\mathbb{Z})$ defined by
\begin{multline}
  \rho(U^{(p,q)})_{a,b}
  \\
  =
  - \mathrm{i} \frac{\sign(q)}{\sqrt{2 \, N \, | q|}} \,
  \mathrm{e}^{
    -\frac{\pi \mathrm{i}}{4} \, \Phi(U^{(p,q)})
  } \,
  \mathrm{e}^{
    \frac{\pi \mathrm{i}}{2 N q} s b^2 
  } \,
  \sum_{
    \substack{
      \gamma \mod 2 N q
      \\
      \gamma = a \mod 2 N}
  }
  \mathrm{e}^{
    \frac{\pi \mathrm{i}}{2 N q}
    p \gamma^2
  } \,
  \left(
    \mathrm{e}^{
      \frac{\pi \mathrm{i}}{ N q}
       \gamma b
    }
    -
    \mathrm{e}^{
      - \frac{\pi \mathrm{i}}{ N q}
      \gamma b
    }
  \right)
\end{multline}
for $1 \leq a, b \leq N-1$~\cite{LCJeff92a}.
This representation is constructed from
\begin{equation}
  \label{affine_su}
  \begin{gathered}
    \rho(S)_{a,b}
    =
    \sqrt{\frac{2}{N}} \,
    \sin
    \left( \frac{ a \, b }{N} \, \pi \right)
    \\[2mm]
    \rho(T)_{a,b}
    =
    \mathrm{e}^{\frac{\pi \mathrm{i}}{2 N} a^2 - \frac{\pi \mathrm{i}}{4}}
    \,
    \delta_{a,b}
  \end{gathered}
\end{equation}
with
\begin{align}
  \label{PSL_ST}
  S
  &=
  \begin{pmatrix}
    0 & -1 \\
    1 & 0
  \end{pmatrix}
  &
  T 
  & =
  \begin{pmatrix}
    1 & 1 \\
    0 & 1
  \end{pmatrix}
\end{align}
satisfying
\begin{equation*}
  S^2 = (S \, T)^3 = 1
\end{equation*}

Based on the fact that the Seifert fibered manifold $\Sigma(\vv{p})$
has a surgery description as in Fig.~\ref{fig:surgery}, we have the
following result.

\begin{prop}[\cite{LawreRozan99a}]
  For the Seifert fibered integral  homology sphere with $M$-exceptional
  fibers
  $\mathcal{M}=\Sigma(\vv{p})$,
  the WRT invariant is given by
  \begin{multline}
    \label{result_Rozansky}
    \E^{\frac{2 \pi \I}{N}
      \left( \frac{\phi(\vv{p})}{4} - \frac{1}{2} \right)
    } \,
    \left(
      \E^{\frac{2 \pi \I}{N}} - 1
    \right) \,
    \tau_N
    (\mathcal{M})
    \\
    =
    \frac{\E^{\frac{\pi \I}{4}}}{
      2 \sqrt{2 \, P \, N}
    } \,
    \sum_{
      \substack{
        n=0 \\
        N \,  \nmid \, n
      }}^{2 P N -1}
        \E^{- \frac{1}{2 P N} n^2 \pi \I } \,
    \frac{
      \prod_{j=1}^M
      \left(
        \E^{\frac{  n}{N p_j} \pi \I}
        -
        \E^{- \frac{ n}{N p_j} \pi \I }
      \right)
    }{
      \left(
        \E^{\frac{n}{N} \pi \I}
        -
        \E^{ - \frac{ n}{N} \pi \I}
      \right)^{M-2}
    }
  \end{multline}
  Here $\phi(\vv{p})$ is defined by
  \begin{equation}
    \label{define_phi_p}
    \phi(\vv{p})
    =
    3-
    \frac{1}{P} + 12 \sum_{j=1}^M
    s\left(
      \frac{P}{p_j} , p_j
    \right)
  \end{equation}
  where $s(b,a)$ denotes the Dedekind sum~\eqref{Dedekind_sum}.
\end{prop}
\begin{proof}[Outline of Proof]
  We use the surgery formula~\eqref{define_invariant}, in which the
  the colored Jones polynomial for a variant of Hopf link
  $\mathcal{L}$ depicted in   Fig.~\ref{fig:surgery} is given by
  \begin{equation*}
    J_{k_0,k_1, \dots, k_M}(\mathcal{L})
    =
    \frac{1}{
      \sin(\pi/N)
    } \,
    \frac{
      \prod_{j=1}^M \sin \left( \frac{k_0 k_j}{N} \pi \right)
    }{
      \left[
        \sin \left( \frac{k_0}{N} \pi \right)
      \right]^{M-1}
    }
  \end{equation*}
  Here $k_0$ is a color of component which has a linking number $1$ with any other components
  of $\mathcal{L}$.
  After some computations using the Gauss sum reciprocity
  formula~\eqref{Gauss_reciprocity}, we get
  \begin{multline}
    \label{WRT_multisum}
    \E^{\frac{2 \pi \I}{N}
      \left( \frac{\phi(\vv{p})}{4} - \frac{1}{2} \right)
    } \,
    \left(
      \E^{\frac{2 \pi \I}{N}} - 1
    \right) \,
    \tau_N
    (\mathcal{M})
    \\
    =
    \frac{\E^{\frac{\pi \I}{4}}}{
      2 \sqrt{2 \, P \, N}
    } \,
    \sum_{k_0 =1}^{N-1} \sum_{n_j \mod p_j}
    \frac{1}{
      \left(
        \E^{\frac{k_0}{N} \pi \I} -
        \E^{-\frac{k_0}{N} \pi \I}
      \right)^{M-2}
    }
    \\
    \times
    \prod_{j=1}^M
    \E^{- \frac{q_j}{p_j} \frac{(k_0+2 N n_j)^2}{2 N}} \,
    \left(
      \E^{\frac{k_0+ 2 N n_j}{N p_j} \pi \I}
      -
      \E^{- \frac{k_0+ 2 N n_j}{N p_j} \pi \I}
    \right)
  \end{multline}
  We see that
  the summand in~\eqref{WRT_multisum} is
  invariant under
  \begin{itemize}
  \item $k_0 \to k_0 + 2 N$ and $n_j \to n_j -1$ for all $j$,
  \item $n_j \to n_j + p_j$ .
  \end{itemize}
  Using that $p_j$ are pairwise coprime,
  we can then rewrite the multi-sum of~\eqref{WRT_multisum}
  into a single-sum of~\eqref{result_Rozansky}.
\end{proof}

The WRT
invariant  can be rewritten in the integral form as follows;
\begin{prop}[\cite{LawreRozan99a}]
  \begin{multline}
    \label{invariant_integral}
    \E^{\frac{2 \pi \I}{N}
      \left( \frac{\phi(\vv{p})}{4} - \frac{1}{2} \right)
    } \,
    \left(
      \E^{\frac{2 \pi \I}{N}} - 1
    \right) \,
    \tau_N
    (\mathcal{M})
    \\
    =
    \frac{\E^{\frac{\pi \I}{4}}}{
      2 \sqrt{2 \, P \, N}
    } \,
    \left(
      - 2 \, \pi \, \I \,
      \sum_{m=0}^{2 P -1}
      \Res_{z=m N}
      \frac{g(z)}{1 - \E^{-2 \pi \I z}}
      +
      \int\limits_{\mathcal{C}} g(z) \, \mathrm{d} z
    \right)
  \end{multline}
  where 
  \begin{equation}
    g(z) = \E^{-\frac{z^2}{2 P N} \pi \I} \,
    \frac{
      \prod_{j=1}^M \left(
        \E^{\frac{z}{ N p_j} \pi \I}
        - \E^{- \frac{z}{N p_j} \pi \I}
      \right)
    }{
      \left(
        \E^{\frac{z}{N} \pi \I} -
        \E^{- \frac{z}{N} \pi \I}
      \right)^{M-2}
    }
  \end{equation}
  and the integration path  $\mathcal{C}$ passes  the origin from
  $(-1+\I)\, \infty$ to $(1-\I) \, \infty$.
\end{prop}
\begin{proof}[Outline of Proof]
  Key identity is
  \begin{equation}
    \Theta_M(x)
    =
    \Theta_M(x-N) \, \E^{4 \pi \I P x}
    + 2 \, \pi \, \I 
    \sum_{m=0}^{2P-1} \Res_{z= m N} h_M(z,x)
    +\sum_{
      \substack{
        n=0\\
        N \nmid n
      }}^{2 P N -1}
    f_M(n,x)
  \end{equation}
  Here the function $\Theta_M(x)$ is defined by
  \begin{equation*}
    \Theta_M(x)
    =
    \int\limits_{\mathcal{C}} h_M(z , x ) \, \mathrm{d} z
  \end{equation*}
  where 
  \begin{equation*}
    h_M(z, x)
    =
    \E^{- \frac{z^2}{2 P N} \pi \I + 2 x \frac{z}{N} \pi \I} \,
    \frac{1}{
      \left(
        \E^{\frac{z}{N} \pi \I}
        -
        \E^{-\frac{z}{N} \pi \I}
      \right)^{M-2}
    } \cdot
    \frac{1}{1 - \E^{-2 \pi \I z}}
    \equiv
    \frac{f_M(z,x)}{
      1 - \E^{-2 \pi \I z}
    }
  \end{equation*}
  As the function $g(z)$ is a linear combination of $f_M(x,z)$, we
  obtain the expression~\eqref{invariant_integral}.
\end{proof}

In view of~\eqref{invariant_integral},
we can decompose the invariant as
\begin{equation}
  \label{tau_res_int}
  \tau_N(\mathcal{M})
  =
  \tau_N^{\mathrm{res}}(\mathcal{M})
  +
  \tau_N^{\mathrm{int}}(\mathcal{M})
\end{equation}
where $  \tau_N^{\mathrm{res}}(\mathcal{M})$
and $\tau_N^{\mathrm{int}}(\mathcal{M})$  respectively
denote  contributions
from the
residue terms and the integral term in~\eqref{invariant_integral}.
It was identified in Ref.~\citen{LawreRozan99a} that 
the residue part $\tau_N^{\mathrm{res}}(\mathcal{M})$
is the contribution from
irreducible flat connections
while the integral term $\tau_N^{\mathrm{int}}(\mathcal{M})$
is the trivial connection contribution.
The trivial connection contribution is related to the Ohtsuki
series~\cite{Ohtsu95b},
and we have as follows;
\begin{prop}
  In the limit $N\to\infty$, the trivial connection contribution
  has the asymptotic expansion as
  \begin{equation}
    \label{Lawrence_find_Ohtsuki}
    \E^{\frac{2 \pi \I}{N} \left( \frac{\phi(\vv{p})}{4} -
        \frac{1}{2}
      \right)
    } \,
    \left(
      \E^{\frac{2 \pi \I}{N}} - 1
    \right) \,
    \tau_N^{\mathrm{int}}(\mathcal{M})
    \simeq
    \sum_{k=0}^\infty \frac{T_{\vv{p}}(k)}{k!} \,
    \left(
      \frac{\pi \, \I}{2 \, P \, N}
    \right)^k
  \end{equation}
  where the $T$-series is given by
  \begin{equation}
    \frac{
      \prod_{j=1}^M
      \sinh
      \left(\frac{P}{p_j } \, x \right)
    }{
      \left[
        \sinh
        \left( P x        \right)
      \right]^{M-2}
    }
    =
    \frac{1}{2}
    \sum_{k=0}^\infty
    \frac{T_{\vv{p}}(k)}{(2 \, k)!} \, x^{2 k}
  \end{equation}
\end{prop}

We will discuss later the relationship with the Ohtsuki series.

Similar  integral with $\tau_N^{\mathrm{int}}(\mathcal{M})$ in the case of the
three exceptional fibers $M=3$ appeared in
studies~\cite{LJMorde33a,GEAndre81e}
of the Ramanujan mock theta functions, which still remain to be mysterious and fascinating
topics.
This  suggests the remarkable fact that the Ramanujan mock theta
functions~\cite{Ramanujan87Book}
are related to the WRT invariant for
the Seifert fibered
manifolds.
See Refs.~\citen{KHikami05a,KHikami05b,KHikami05c} for detail.

The dominating term of the WRT invariant in a limit $N\to \infty$
follows from the irreducible flat connection contributions
$\tau_N^{\mathrm{res}}(\mathcal{M})$ as was expected from the saddle
point approximation~\eqref{expected_asymptotics}.
In terms of the Witten partition function $Z_k(\mathcal{M})$, the
asymptotic expansion of the invariant is given by
\begin{equation*}
  Z_{k-2}(\mathcal{M})
  \simeq
  \sum_{a=0}^{M-3} N^{M-3-a} \, Z_{N-2}^{(a)}(\mathcal{M})
  + \text{trivial connection contribution}
\end{equation*}
and the dominating term $Z_{N-2}^{(0)}(\mathcal{M})$ in $N\to\infty$
can be  computed
as follows when we use an identity
\begin{equation*}
  \sin \left(\pi \, z \right)
  =
  \pi \, z \, \prod_{n=1}^\infty \left( 1 - \frac{z^2}{n^2}\right)
\end{equation*}

\begin{prop}
  In the large $N$ limit, 
  the Witten  partition function $Z_N(\mathcal{M})$ for the
  $M$-exceptional fibered Seifert homology sphere
  $\mathcal{M}=\Sigma(\vv{p})$ is
  dominated by
  \begin{multline}
    \label{asymptotics_Lawrence}
    Z_{N-2}(\mathcal{M})
    \sim
    N^{M-3} \, \frac{2^{M-2}}{(M-2)! \, \sqrt{P}} \,
    \E^{- \frac{\phi(\vv{p})}{2 N} \pi \I} \,
    \E^{-\frac{2 M-3}{4}\pi \I}
    \\
    \times
    \sum_{m=0}^{2P-1} (-1)^{m M} 
    B_{M-2} \left( \frac{m}{2 \, P}\right) \,
    \E^{- \frac{m^2}{2 P} \pi \I N}
    \,
    \left[
      \prod_{j=1}^M
      \sin
      \left( \frac{m}{p_j} \, \pi \right)
    \right]
  \end{multline}
  where $B_k(x)$ is the $k$-th Bernoulli
  polynomial~\eqref{generate_Bernoulli}.
\end{prop}

Among a sum of $2 \, P$ terms  in the right hand side
of~\eqref{asymptotics_Lawrence}, we can classify the summation
by the Chern--Simons invariant, which corresponds to an exponential factor
$-\frac{m^2}{4 \, P} \mod 1$ as will be discussed
later.


\section{Modular Forms and Eichler Integral}
\label{sec:modular_form}

We introduce the vector modular forms with half-integral weight which
play a crucial role in analysis of the WRT invariants for the Seifert
fibered homology spheres.

\subsection{Vector Modular Forms with Half-Integral Weight}
We set $M$-tuple
\begin{equation}
  \vv{\ell}=
  (\ell_1, \dots, \ell_M)
\end{equation}
where $\ell_j$ are integers satisfying
$0 < \ell_j < p_j$.
As before we assume that $p_j$ are pairwise coprime positive integers.
For $M$-tuple $\vv{\ell}$,
we define the periodic function
$\chi_{2 P}^{\vv{\ell}}(n)$ with modulus $2 \, P$ by
\begin{equation}
  \label{define_chi}
  \chi_{2 P}^{\vv{\ell}}(n)
  =
  \begin{cases}
    \displaystyle
    - \prod_{j=1}^M \varepsilon_j ,
    &
    \text{
      if $\displaystyle
      n = P \, \left(
        1 + \sum_{j=1}^M \varepsilon_j \,  \frac{\ell_j}{p_j}
      \right) \mod 2 \, P$}
    \\[4mm]
    0 , &
    \text{otherwise}
  \end{cases}
\end{equation}
where $\varepsilon_j = \pm 1$ for $\forall j$.
We see  that
$\chi_{2 P}^{\vv{\ell}}(n)$ is even (resp. odd)
when $M$ is even (resp. odd),
\begin{equation}
  \chi_{2P}^{\vv{\ell}}(-n)
  = (-1)^M \, \chi_{2 P}^{\vv{\ell}}(n)
\end{equation}
and that it has a mean value zero,
\begin{equation}
  \label{chi_mean_zero}
  \sum_{n=0}^{2 P-1} \chi_{2 P}^{\vv{\ell}}(n) = 0
\end{equation}

We define an involution $\sigma_j$ on $M$-tuple $\vv{\ell}$ by
\begin{equation}
  \sigma_j(\vv{\ell})
  =
  (\ell_1, \dots,\ell_{j-1} ,  p_j - \ell_j , \ell_{j+1} , \dots, \ell_M)
\end{equation}
for $1 \leq j \leq M$.
As we have
\begin{equation}
  \label{sigma2_chi}
  \chi_{2 P}^{\sigma_i \sigma_j(\vv{\ell})}(n)
  =
  \chi_{2 P}^{\vv{\ell}}(n)
\end{equation}
for $1 \leq i, j \leq M$,
the number
of the independent periodic functions
$\chi_{2 P}^{\vv{\ell}}(n)$ for given $\vv{p}$
is
\begin{equation}
  \label{define_Dp}
  D=
  D(\vv{p})
  =
  \frac{1}{2^{M-1}} \,
  \prod_{j=1}^M
  \left(
    p_j -1
  \right)
\end{equation}
We note that
\begin{equation}
  \chi_{2 P}^{\vv{\ell}}(n+P)
  =
  -
  \chi_{2 P}^{\sigma_j(\vv{\ell})}(n)
\end{equation}
for $1 \leq j \leq M$.

We set
\begin{equation}
  q = \exp \left( 2 \, \pi \, \I \, \tau \right)
\end{equation}
where $\tau$ is in the upper half plane,
$\tau \in \mathbb{H}$.
By use of the periodic functions~\eqref{define_chi}
we define  the $q$-series by
\begin{equation}
  \Phi_{\vv{p}}^{\vv{\ell}}(\tau)
  =
  \frac{1}{2} \,
  \sum_{n \in \mathbb{Z}}
  n^{\eveodd(M)} \,
  \chi_{2 P}^{\vv{\ell}}(n)
  \,
  q^{\frac{n^2}{4 P}}
\end{equation}
where we mean
\begin{equation}
  \eveodd(M)
  =
  \frac{1-(-1)^M}{2}
  =
  M \mod 2
  =
  \begin{cases}
    0 & \text{when $M$ is even}
    \\[2mm]
    1 & \text{when $M$ is odd}
  \end{cases}
\end{equation}

The $q$-series $\Phi_{\vv{p}}^{\vv{\ell}}(\tau)$ is proved to be a
vector modular form with half-integral weight.
The $T$-transformation is trivial, and
by use of the Poisson summation formula,
\begin{equation}
  \sum_{n \in \mathbb{Z}} f(n)
  =
  \sum_{n \in \mathbb{Z}} \int_{-\infty}^\infty \E^{-2 \pi \I t n} \,
  f(t) \,
  \mathrm{d} t
\end{equation}
we obtain the transformation formula under the
$S$-transformation as follows.

\begin{prop}
  The $q$-series
  $\Phi_{\vv{p}}^{\vv{\ell}}(\tau)$
  is a vector modular form with
  weight
  $3/2$ (resp. $1/2$)
  when $M$ is even (resp. odd).
  Under the $S$- and $T$-transformations~\eqref{PSL_ST} we have
  \begin{gather}
    \label{Phi_under_S}
    \Phi_{\vv{p}}^{\vv{\ell}}(\tau)
    =
    \left(
      \frac{  \  \I  \  }{\tau}
    \right)^{\frac{3}{2} - \eveodd(M)} \,
    \sum_{
      \vv{\ell^\prime}
    }
    \mathbf{S}_{\vv{\ell^\prime}}^{\vv{\ell}}
    \,
    \Phi_{\vv{p}}^{\vv{\ell^\prime}}
    (-1/\tau)
    \\[2mm]
    \Phi_{\vv{p}}^{\vv{\ell}}(\tau+1)
    =
    \mathbf{T}^{\vv{\ell}}
    \,
    \Phi_{\vv{p}}^{\vv{\ell}}(\tau)
  \end{gather}
  Here a sum of $\vv{\ell^\prime}$ runs over
  $D$-dimensional
  space~\eqref{define_Dp}, and matrix elements of
  $D\times D$ matrices
  $\mathbf{S}$ and $\mathbf{T}$ are
  respectively given by
  \begin{gather}
    \label{S_matrix}
    \mathbf{S}_{\vv{\ell^\prime}}^{\vv{\ell}}
    =
    \frac{2^M \,
      \I^{M-\eveodd(M)}
    }{\sqrt{2 \, P}} 
    \,
    (-1)^{
      P
      \left(
        1 + \sum_{j=1}^M \frac{\ell_j+\ell_j^\prime}{p_j}
      \right)
      +
      P \sum_j \sum_{k \neq j }
      \frac{\ell_j \, \ell_k^\prime}{
        p_j  \, p_k
      }
    } \,
    \prod_{j=1}^M
    \sin
    \left(
      P \, \frac{\ell_j \, \ell_j^\prime}{p_j^{~2}} \,\pi
    \right)
    \\[2mm]
    \mathbf{T}^{\vv{\ell}}
    =
    \exp
    \left(
      \frac{P}{2} \,
      \Bigl(
      1+ \sum_{j=1}^M \frac{\ell_j}{p_j}
      \Bigr)^2 \,
      \pi \, \I
    \right)
    \label{T_matrix}
  \end{gather}
\end{prop}

Our  vector modular form $\Phi_{\vv{p}}^{\vv{\ell}}(\tau)$
may be a generalization of modular forms which appear 
as the character of the affine Lie algebra
$\widehat{su(2)}$ (a case of $M=1$, and $\Psi_{P}^{(a)}(\tau)$ defined below)
and  as the minimal Virasoro model (a case of $M=2$)
up to the power of the Dedekind $\eta$ function.

For our later use,
we introduce other families of vector modular forms.
We define the even periodic functions by
\begin{gather}
  \theta_{2 P}^{(a)}(n)
  =
  \begin{cases}
    1 & \text{for $n = \pm a \mod 2 \, P$}
    \\[2mm]
    0 & \text{otherwise}
  \end{cases}
\end{gather}
for $0 \leq a \leq P$, and the odd periodic function by
\begin{gather}
  \psi_{2 P}^{(a)}(n)
  =
  \begin{cases}
    \pm 1 &
    \text{for $n \equiv \pm a \mod 2 \, P$}
    \\[2mm]
    0 &
    \text{otherwise}
  \end{cases}
\end{gather}
for
$ 0 < a < P$.
We then define two families of $q$-series by
\begin{gather}
  \Theta_P^{(a)}(\tau)
  =
  \frac{1}{2} \sum_{n \in \mathbb{Z}}
  \theta_{2 P}^{(a)}(n) \,
  q^{\frac{n^2}{4 P}}
  \\[2mm]
  \Psi_P^{(a)}(\tau)
  =
  \frac{1}{2}
  \sum_{n \in \mathbb{Z}}
  n \, \psi_{2 P}^{(a)}  (n) \,
  q^{\frac{n^2}{4 \, P}}
\end{gather}
These families are also vector  modular forms with half-integral weight.
Namely
we see that $\Theta_P^{(a)}(\tau)$ is a vector modular form with
weight $1/2$ satisfying
\begin{gather}
  \Theta_P^{(a)}(\tau)
  =
  \sqrt{
    \frac{  \  \I  \  }{\tau}
  } \,
  \sum_{b=0}^{P}
  \mathbf{N}_b^a \,
  \Theta_P^{(b)}(-1/\tau)
  \\[2mm]
  \Theta_P^{(a)}(\tau+1)
  =
  \exp
  \left(
    \frac{a^2}{2 \, P} \, \pi \, \I
  \right) \,
  \Theta_P^{(a)}(\tau)
\end{gather}
where $\mathbf{N}$ is $(P+1)\times(P+1)$ matrix defined by
\begin{equation}
  \mathbf{N}^a_b
  =
  \begin{cases}
    \displaystyle
    \frac{1}{\sqrt{2 \, P}} &
    \text{for $a=0$}
    \\[4mm]
    \displaystyle
    \sqrt{\frac{2}{P}} 
    \cos \left(
      \frac{a \, b}{P} \, \pi
    \right)
    & \text{for $a \neq 0, P$}
    \\[4mm]
    \displaystyle
    \frac{1}{\sqrt{2 P}}
    \cos\left( b \, \pi \right)
    & \text{for $a = P$}
  \end{cases}
\end{equation}
The vector modular form $\Psi_P^{(a)}(\tau)$ with weight $3/2$
fulfills
the following transformation formulae;
\begin{gather}
  \label{Psi_under_S}
  \Psi_P^{(a)}(\tau)
  =
  \left(
    \frac{  \  \I  \  }{\tau}
  \right)^{3/2} \,
  \sum_{b=1}^{P-1}
  \mathbf{M}_b^a \,
  \Psi_P^{(b)}(-1/\tau)
  \\[2mm]
  \Psi_P^{(a)}(\tau+1)
  =
  \exp
  \left(
    \frac{a^2}{2 \, P} \, \pi \, \I
  \right) \,
  \Psi_P^{(a)}(\tau)
\end{gather}
where $\mathbf{M}$ is a $(P-1)\times (P-1)$ matrix defined by
\begin{equation}
  \mathbf{M}_b^a =
  \sqrt{\frac{  \  2  \  }{P}} \,
  \sin
  \left(\frac{a \, b}{P} \, \pi \right)
\end{equation}

\subsection{Eichler Integrals}

The Eichler integral is originally defined for  modular forms with
integral weight $\geq 2$
(see, \emph{e.g.}, Ref.~\citen{SLang76Book}).
In our cases, the vector modular forms
$\Phi_{\vv{p}}^{\vv{\ell}}(\tau)$ have a half-integral weight,
so we follow a method of Refs.~\citen{DZagie01a,LawrZagi99a}
to  define a variant of the Eichler integrals.

We define the Eichler integrals
$\widetilde{\Phi}_{\vv{p}}^{\vv{\ell}}(\tau)$
of the vector modular form $\Phi_{\vv{p}}^{\vv{\ell}}(\tau)$
by
\begin{equation}
  \label{define_Eichler}
  \widetilde{\Phi}_{\vv{p}}^{\vv{\ell}}(\tau)
  =
  \sum_{n=0}^\infty
  n^{1-\eveodd(M)} \,
  \chi_{2 P}^{\vv{\ell}}(n) \,
  q^{\frac{n^2}{4 P}}
\end{equation}
This can be regarded as a ``half-derivative'' (resp. ``half-integral'') of
the modular form $\Phi_{\vv{p}}^{\vv{\ell}}(\tau)$ with respect to
$\tau$ when $M$ is even (resp. odd).
When $M$ is odd, the $q$-series
$\widetilde{\Phi}_{\vv{p}}^{\vv{\ell}}(\tau)$ might be called the
false theta function {\`{a}} la Rogers~\cite{LJRogers17a}.

\begin{prop}
  We assume $N_1 $ and $N_2$ are coprime integers, and $N_1>0$.
  Limiting values of the Eichler integrals
  $\widetilde{\Phi}_{\vv{p}}^{\vv{\ell}}(\tau)$
  in $\tau \to N_2/N_1$ are given as follows;
  \begin{equation}
    \label{limiting_value_Eichler}
    \widetilde{\Phi}_{\vv{p}}^{\vv{\ell}}(
    N_2 / N_1
    )
    =
    -
    \left( P \, N_1 \right)^{1-\eveodd(M)} 
    \sum_{k=0}^{2 P N_1}
    \chi_{2 P}^{\vv{\ell}}(k) \,
    \E^{
      \frac{N_2}{N_1} \, \frac{k^2}{2 P} \,
      \pi \,
      \I
    } \,
    B_{2-\eveodd(M)}
    \left(
      \frac{k}{2 \, P \, N_1}
    \right)
  \end{equation}
  where
  $B_k(x)$ is the $k$-th Bernoulli polynomial.
%
\end{prop}

We note that
\begin{equation*}
  \widetilde{\Phi}_{\vv{p}}^{\vv{\ell}}(\tau+1)
  =
  \mathbf{T}^{\vv{\ell}} \,
  \widetilde{\Phi}_{\vv{p}}^{\vv{\ell}}(\tau)
\end{equation*}
from which we have for $N \in \mathbb{Z}$
\begin{align*}
  \widetilde{\Phi}_{\vv{p}}^{\vv{\ell}}(N)
  & =
   \left( \mathbf{T}^{\vv{\ell}} \right)^N \, \widetilde{\Phi}_{\vv{p}}^{\vv{\ell}}(0)
  \\
  & =
  - \left( \mathbf{T}^{\vv{\ell}} \right)^N \,
  P^{1-\eveodd(M)} \,
  \sum_{k=0}^{2P}
  \chi_{2P}^{\vv{\ell}}(k) \,
  B_{2-\eveodd(M)}
  \left(\frac{k}{2 \, P} \right)
\end{align*}

A proof follows straightforwardly when we use the following
lemma~\cite{LawrZagi99a}
(see also Ref.~\citen{KOno04Book}).
\begin{lemma}
  \label{lemma:classical_limit}
  Let $C_f(n)$ is a periodic function with modulus $f$
  and mean  value zero.
  Then we have as $t \searrow 0$
  \begin{gather}
    \sum_{n=1}^\infty
    C_f(n) \, \E^{-n^2 t}
    \simeq
    \sum_{n=0}^\infty
    L\left( -2 \, n, C_f\right) \,
    \frac{(-t)^n}{n!}
    \\[2mm]
    \sum_{n=1}^\infty
    n \, C_f(n) \, \E^{-n^2 t}
    \simeq
    \sum_{n=0}^\infty
    L\left( -2 \, n - 1, C_f\right) \,
    \frac{(-t)^n}{n!}
  \end{gather}
  where the $L$-function is
  \begin{align}
    L\left(s , C_f \right)
    & =
    \sum_{n=1}^\infty \frac{C_f(n)}{n^s}
    \\
    & = f^{-s} \sum_{k=1}^f C_f(k) \, 
    \zeta \left(s, \frac{k}{f} \right)
    \nonumber
  \end{align}
\end{lemma}
Note that  the Hurwitz zeta function $\zeta(s,z)$ defined by
\begin{equation}
  \zeta(s,z) = \sum_{n=0}^\infty
  \frac{1}{(n+z)^s}
\end{equation}
has
an analytic continuation 
for $k \in \mathbb{Z}_{>0}$ as
\begin{equation}
  \label{Hurwitz_Bernoulli}
  \zeta(1-k, z)
  = -
  \frac{B_k(z)}{k}
\end{equation}

\begin{prop}
  The limiting values
  $\widetilde{\Phi}_{\vv{p}}^{\vv{\ell}}(\alpha)$ 
  of the Eichler integrals
  with 
  $\alpha\in  \mathbb{Q}$
  satisfy a nearly
  modular property.
  In a  limit $N\to\infty$, we have the transformation
  formula as an asymptotic expansion  as follows;
  \begin{multline}
    \label{nearly_modular_Eichler}
    \widetilde{\Phi}_{\vv{p}}^{\vv{\ell}}(1/N)
    +
    \left(
      \frac{N}{\I}
    \right)^{\frac{3}{2} - \eveodd(M)} \,
    \sum_{\vv{\ell^\prime}}
    \mathbf{S}_{\vv{\ell^\prime}}^{\vv{\ell}}
    \,
    \widetilde{\Phi}_{\vv{p}}^{\vv{\ell^\prime}}(-N)
    \\
    \simeq
    \sum_{k=0}^\infty
    \frac{
      L \left(
        -2 \, k - 1 + \eveodd(M) , \chi_{2 P}^{\vv{\ell}}
      \right)
    }{
      k!
    } \,
    \left(
      \frac{
        \pi \, \I
      }{
        2 \, P \, N
      }
    \right)^k
  \end{multline}
%
  Here $N \in \mathbb{Z}$,
  and
  a sum of $M$-tuples $\vv{\ell^\prime}$
  runs over
  $D$-dimensional space.
\end{prop}
\begin{proof}
  We introduce another variant of the Eichler integral by
  \begin{equation}
    \widehat{\Phi}_{\vv{p}}^{\vv{\ell}}(z)
    =
    \begin{cases}
      \displaystyle
      -\sqrt{\frac{P \, \I}{2 \, \pi^2}} \,
      \int_{\bar{z}}^\infty
      \frac{
        \Phi_{\vv{p}}^{\vv{\ell}}(\tau)
      }{
        \left( \tau - z \right)^{3/2}
      } \,
      \mathrm{d} \tau ,
      &
      \text{when $M$ is even}
      \\
      \\
      \displaystyle
      \frac{1}{\sqrt{2 \, P \, \I}} \,
      \int_{\bar{z}}^\infty
      \frac{
        \Phi_{\vv{p}}^{\vv{\ell}}(\tau)
      }{
        \sqrt{\tau - z}
      }
      \,
      \mathrm{d} \tau ,
      &
      \text{when $M$ is odd}
    \end{cases}
  \end{equation}
  where we assume that $z$ is in the lower half-plane,  $z \in \mathbb{H}^-$,
  and $\bar{z}$ denotes a complex conjugate of $z$.
  By use of the $S$-transformation~\eqref{Phi_under_S}, we have
  \begin{equation}
    \label{nearly_modular_widehat}
    \widehat{\Phi}_{\vv{p}}^{\vv{\ell}}(z)
    +
    \left(
      \frac{1}{
        \I \, z   
      }
    \right)^{\frac{3}{2} - \eveodd(M)}
    \,
    \sum_{\vv{\ell^\prime}}
    \mathbf{S}_{\vv{\ell^\prime}}^{\vv{\ell}}
    \,
    \widehat{\Phi}_{\vv{p}}^{\vv{\ell^\prime}}(-1/z)
    =
    r_{\Phi_{\vv{p}}^{\vv{\ell}}}(z; 0)
  \end{equation}
  Here $r_{\Phi_{\vv{p}}^{\vv{\ell}}}(z;\alpha)$
  is an analogue of the period function defined by
  \begin{equation}
    r_{\Phi_{\vv{p}}^{\vv{\ell}}}(z; \alpha)
    =
    \begin{cases}
      \displaystyle
      -\sqrt{\frac{P \, \I}{2 \, \pi^2}} \,
      \int_\alpha^\infty
      \frac{
        \Phi_{\vv{p}}^{\vv{\ell}}(\tau)
      }{
        \left( \tau - z \right)^{3/2}
      } \,
      \mathrm{d} \tau  ,
      &
      \text{when $M$ is even}
      \\
      \\
      \displaystyle
      \frac{1}{
        \sqrt{2 \, P \, \I}
      } \,
      \int_\alpha^\infty
      \frac{\Phi_{\vv{p}}^{\vv{\ell}}(\tau)}{
        \sqrt{\tau-z}
      } \,
      \mathrm{d} \tau ,
      &
      \text{when $M$ is odd}
    \end{cases}
  \end{equation}
  where $\alpha \in \mathbb{Q}$.
We find that $\widehat{\Phi}_{\vv{p}}^{\vv{\ell}}(z)$  takes the same
limiting value with that of $\widetilde{\Phi}_{\vv{p}}^{\vv{\ell}}(\tau)$ in a
limit $z, \tau \to \alpha \in \mathbb{Q}$;
\begin{equation}
  \label{tilde_and_hat}
  \left.
    \widetilde{\Phi}_{\vv{p}}^{\vv{\ell}}(\tau)
  \right|_{\tau \to \alpha}
  =
  \left.
    \widehat{\Phi}_{\vv{p}}^{\vv{\ell}}(z)
  \right|_{z \to \alpha}
\end{equation}
The right hand side of~\eqref{nearly_modular_Eichler} arises from
an asymptotic expansion of $r_{\Phi_{\vv{p}}^{\vv{\ell}}}(z; 0)$,
and we obtain~\eqref{nearly_modular_Eichler}.
\end{proof}

For our later use, we study differentials of the Eichler
integral $\widetilde{\Phi}_{\vv{p}}^{\vv{\ell}}(\tau)$,
\emph{i.e.}, ``fractional derivatives''
of the vector modular form $\Phi_{\vv{p}}^{\vv{\ell}}(\tau)$.
From the definition~\eqref{define_Eichler} of the Eichler integral, we
have
for $b \in \mathbb{Z}_{\geq 0}$
\begin{equation}
  \left(
    \frac{2 \, P}{\pi \, \I} \,
    \frac{\mathrm{d}}{\mathrm{d} \tau}
  \right)^b \,
  \widetilde{\Phi}_{\vv{p}}^{\vv{\ell}}
  (\tau)
  =
  \sum_{n=0}^\infty
  n^{2b+1 - \eveodd(M)} \, \chi_{2 P}^{\vv{\ell}}(n) \,
  q^{\frac{n^2}{4 P}}
\end{equation}
By the same computation with~\eqref{limiting_value_Eichler}, we have
the following.
\begin{prop}
  The limiting values of \emph{fractional} derivative of the vector
  modular forms $\Phi_{\vv{p}}^{\vv{\ell}}(\tau)$ are given by
  \begin{multline}
    \label{limiting_frac_differential}
    \left.
      \left(
        \frac{2 \, P}{\pi \, \I} \,
        \frac{\mathrm{d}}{\mathrm{d} \tau}
      \right)^b \,
      \widetilde{\Phi}_{\vv{p}}^{\vv{\ell}}
      (\tau)
    \right|_{\tau\to \frac{N_2}{N_1}}
    \\
    =
    -
    \frac{\left( 2 \, P \, N_1 \right)^{2b+1 -\eveodd(M)}}{
      2\, b + 2 - \eveodd(M) }
    \sum_{n=1}^{2 P N_1}
    \chi_{2 P}^{\vv{\ell}}(n) \,
    \E^{\frac{N_2}{N_1} \frac{n^2}{2 P} \pi \I}
    \,
    B_{2 b +2 - \eveodd(M)}
    \left(
      \frac{n}{2 \, P \, N_1}
    \right)
  \end{multline}
  where $N_1>0 $ and $N_2$ are coprime integers.
\end{prop}

The nearly modular property~\eqref{nearly_modular_Eichler} of the
Eichler integral
gives the following asymptotic expansion of~\eqref{limiting_frac_differential}.
\begin{prop}
  In the limit $N\to\infty$, we have the
  following asymptotic expansion;
  \begin{multline}
    \label{modular_frac_differential}
    N^{2 b +1 - \eveodd(M)}
    \sum_{n=1}^{2 P N}
    \chi_{2 P}^{\vv{\ell}}(n) \,
    \E^{\frac{n^2}{2 P N} \pi \I} \,
    B_{2 b + 2 - \eveodd(M)}
    \left(
      \frac{n}{2 \, P\, N}
    \right)
    \\
    \simeq
    \frac{-1}{\I^{\frac{3}{2} - \eveodd(M)}}
    \sum_{j=0}^b
    N^{b+j+\frac{3}{2} - \eveodd(M)} \,
    \left(
      \frac{\I}{
        2 \, P \, \pi 
      }
    \right)^{b-j} \,
    K_{b,\eveodd(M)}^{(j)}
    \\
    \times
    \frac{2 \, b + 2 - \eveodd(M)}{
      2 \, j + 2 - \eveodd(M)}
    \sum_{\vv{\ell^\prime}}
    \mathbf{S}_{\vv{\ell^\prime}}^{\vv{\ell}}
    \left[
      \sum_{n=1}^{2 P}
      \chi_{2 P}^{\vv{\ell^\prime}}(n) \,
      B_{2 j + 2 - \eveodd(M)}
      \left(
        \frac{n}{2 \, P}
      \right) 
\right] \,
      \E^{- 
        \frac{P}{2} 
        \left(
        1+\sum_j \frac{\ell_j^\prime}{p_j}
        \right)^2 \pi \I N
      }
    \\
    -
    \frac{2 \, b+ 2 - \eveodd(M)}{
      \left(2 \, P\right)^{
        2 b + 1 - \eveodd(M)
      }}
    \sum_{k=0}^\infty
    \frac{
      L \left(
        -2 \, k - 2 \, b - 1 + \eveodd(M) , \chi_{2 P}^{\vv{\ell}}
      \right)
    }{
      k!}
    \,
    \left(
      \frac{\pi \, \I}{2 \, P \, N}
    \right)^k
  \end{multline}
  where the sum of $M$-tuples $\vv{\ell^\prime}$ runs over
  $D$-dimensional space, and we have
  \begin{equation}
    K_{b,x}^{(j)}
    =
    \begin{pmatrix}
      b \\
      j
    \end{pmatrix}
    \,
    \prod_{k=0}^{b-j-1}
    \left(
      \frac{1}{2} + b - x - k
    \right)
    \label{define_K-number}
  \end{equation}
\end{prop}
\begin{proof}
  We differentiate~\eqref{nearly_modular_widehat} with respect to
  $z$, and then take a limit $z\to 1/N$.
  We get
  \begin{multline}
    \label{differ_before_limit}
    \left.
      \frac{\mathrm{d}^b}{\mathrm{d} z^b}
      \widehat{\Phi}_{\vv{p}}^{\vv{\ell}}(z)
    \right|_{z \to \frac{1}{N}}
    +
    \frac{1}{
      \I^{\frac{1}{2} + \eveodd(M)}
    }
    \sum_{\vv{\ell^\prime}}
    \mathbf{S}_{\vv{\ell^\prime}}^{\vv{\ell}} \,
    \left.
      \left(
        w^2 \, \frac{\mathrm{d}}{\mathrm{d} w}
      \right)^b \,
      w^{\frac{3}{2} - \eveodd(M)} \,
      \widehat{\Phi}_{\vv{p}}^{\vv{\ell^\prime}}(w)
    \right|_{w \to -N}
    \\
    \simeq
    \left(
      \frac{\pi  \, \I}{2 \, P}
    \right)^b
    \sum_{k=0}^\infty
    \frac{
      L \left(
        - 2 \, k - 2 \, b -1 + \eveodd(M), \chi_{2 P}^{\vv{\ell}}
      \right)
    }{
      k!
    }  \,
    \left(
      \frac{\pi \, \I}{2 \, P \, N}
    \right)^k
  \end{multline}
  In the left hand side, we use
  \begin{equation*}
    \left(
      w^2 \, \frac{\mathrm{d}}{\mathrm{d} w}
    \right)^n
    =
    \sum_{m=1}^n
    A_n^{(m)} \, w^{n+m} \,
    \frac{\mathrm{d}^m}{\mathrm{d} w^m}
  \end{equation*}
  where
  $A_n^{(m)}$  is the Lah number
  defined by
  \begin{equation*}
    A_n^{(m)}
    =
    \frac{n!}{m!} \,
    \begin{pmatrix}
      n-1 \\
      m-1
    \end{pmatrix}
  \end{equation*}
  This denotes the number of partitions of $\{1,2,\dots,n\}$
  into $m$ \emph{lists}
  (a ``list'' denotes an ordered subset)~\cite{NSloanWWW,Stanley97Book},
  and satisfies the recursion relation
  \begin{equation*}
    A_{n+1}^{(m)}
    = A_n^{(m-1)} + ( n+ m) \, A_n^{(m)}
  \end{equation*}
  Thus we have
  \begin{equation*}
    \left( w^2 \, \frac{\mathrm{d}}{\mathrm{d} w} \right)^b \,
    w^{\frac{3}{2} - \eveodd(M)} \,
    \widehat{\Phi}_{\vv{p}}^{\vv{\ell}}(w)
    =
    \sum_{j=0}^b w^{\frac{3}{2} - \eveodd(M)+b+j} \,
    K_{b,\eveodd(M)}^{(j)} \,
    \frac{\mathrm{d}^j \widehat{\Phi}_{\vv{p}}^{\vv{\ell}} (w) }{\mathrm{d} w^j}
  \end{equation*}
  Here the $K$-number is computed as
  \begin{equation*}
    K_{b,x}^{(j)}
    =
    \sum_{k=j}^b A_b^{(k)} \,
    \begin{pmatrix}
      k \\
      j
    \end{pmatrix} \,
    \frac{
      \Gamma \left(\frac{5}{2} - x \right)
    }{
      \Gamma \left( \frac{5}{2} -x - k+j \right)
    }
  \end{equation*}
  which reduces to~\eqref{define_K-number}
  applying~\eqref{generate_Stirling} and~\eqref{Stirling_sum_sum}.
  
  As we have~\eqref{tilde_and_hat}, we
  get~\eqref{modular_frac_differential} with a help
  of~\eqref{limiting_frac_differential}.
\end{proof}

In the same method, we have formulae of the asymptotic expansions
concerning to the periodic functions with mean value zero.

\begin{coro}
  \label{coro:asymptotic_Bernoulli}
  We assume that $C_{2P}(n)$ is an odd or even
  periodic function with modulus
  $2 \,  P$, and that it
  satisfies
  \begin{itemize}
  \item 
    a mean value zero condition,
    \begin{equation*}
      \sum_{n=0}^{2P-1} C_{2P}(n) = 0
    \end{equation*}

  \item
    $C_{2P}(0)=0$.
  \end{itemize}
  In the limit $N\to\infty$, we have the following asymptotic
  expansions;
  \begin{itemize}
  \item   $C_{2P}(n)$ is odd;
    \begin{multline}
      \label{asymptotic_psi_Bernoulli}
      N^{2 b}
      \sum_{n=1}^{2 P N}
      C_{2 P}(n) \,
      \E^{\frac{n^2}{2 P N} \pi \I} \,
      B_{2 b + 1}
      \left(
        \frac{n}{2 \, P\, N}
      \right)
      \\
      \simeq
      \frac{-2}{\I^{\frac{1}{2} }}
      \sum_{j=0}^b
      N^{b+j+\frac{1}{2}} \,
      \left(
        \frac{\I}{
          2 \, P \, \pi 
        }
      \right)^{b-j} \,
      K_{b,  1}^{(j)}      \,
      \frac{2 \, b + 1}{
        2 \, j + 1}
      \sum_{a=1}^{P-1}      C_{2P}(a) \,
      \sum_{c=1}^{P-1}
      \mathbf{M}^{a}_{c}
      B_{2 j + 1}
      \left(
        \frac{c}{2 \, P}
      \right) \,
      \E^{-N 
        \frac{c^2}{2 P}
        \pi \I
      }
      \\
      -
      \frac{2 \, b+ 1}{
        \left(2 \, P\right)^{
          2 b 
        }}
      \sum_{k=0}^\infty
      \frac{
        L \left(
          -2 \, k - 2 \, b  , C_{2 P}
        \right)
      }{
        k!}
      \,
      \left(
        \frac{\pi \, \I}{2 \, P \, N}
      \right)^k
    \end{multline}

  \item   $C_{2P}(n)$ is even;
    \begin{multline}
      \label{asymptotic_theta_Bernoulli}
      N^{2 b +1 }
      \sum_{n=1}^{2 P N}
      C_{2 P}(n) \,
      \E^{\frac{n^2}{2 P N} \pi \I} \,
      B_{2 b + 2 }
      \left(
        \frac{n}{2 \, P\, N}
      \right)
      \\
      \simeq
      \frac{-1}{\I^{\frac{3}{2} }}
      \sum_{j=0}^b
      N^{b+j+\frac{3}{2}} \,
      \left(
        \frac{\I}{
          2 \, P \, \pi 
        }
      \right)^{b-j} \,
      K_{b, 0}^{(j)}
      \,
      \frac{ b + 1}{   j + 1}
      \sum_{a=1}^P C_{2P}(a)
      \\
      \times
      \sum_{c=0}^P
      \mathbf{N}^{a}_{c} \,
      \left( 2 - \delta_{c,0} - \delta_{c,P} \right) \,
      B_{2 j + 2 }
      \left(
        \frac{c}{2 \, P}
      \right)
      \,
      \E^{-N 
        \frac{c^2}{2 P} \pi \I
      }
      \\
      -
      \frac{2 \, b+ 2}{
        \left(2 \, P\right)^{
          2 b + 1 
        }}
      \sum_{k=0}^\infty
      \frac{
        L \left(
          -2 \, k - 2 \, b - 1  , C_{2 P}
        \right)
      }{
        k!}
      \,
      \left(
        \frac{\pi \, \I}{2 \, P \, N}
      \right)^k
    \end{multline}
  \end{itemize}
\end{coro}
\begin{proof}
  We use that 
  \begin{equation*}
    C_{2P}(n) =
    \begin{cases}
      \displaystyle
      \sum_{a=1}^{P-1} C_{2P}(a) \, \psi_{2P}^{(a)}(n) 
      & \text{when $C_{2P}(n)$ is odd}
      \\[4.8mm]
      \displaystyle
      \sum_{a=1}^P C_{2P}(a) \, \theta_{2P}^{(a)}(n)      
      & \text{when $C_{2P}(n)$ is even}
    \end{cases}
  \end{equation*}
  In both cases,  due to a condition $C_{2P}(0)=0$, we see
  that periodic functions which appear in the modular transformation
  formula such as~\eqref{differ_before_limit}
  have a mean value zero.
  Applying Lemma~\ref{lemma:classical_limit}, we obtain asymptotic
  expansions.
\end{proof}

\section{Exact Asymptotic Expansion of the WRT invariants}
\label{sec:asymptotic}

\subsection{Exact Asymptotic Expansion}

As a preparation to obtain the exact asymptotic expansion of the WRT
invariant $\tau_N(\mathcal{M})$ for the $M$-exceptional fibered
Seifert homology
sphere $\mathcal{M}=\Sigma(\vv{p})$ by use of the vector modular form
$\Phi_{\vv{p}}^{\vv{\ell}}(\tau)$
defined in previous section,
we give $q$-series identities at the root of unity
related to  the Bernoulli polynomials.

\begin{lemma}
  We set $\omega_N$ as the $N$-th primitive root of unity;
  \begin{equation*}
    \omega_N=\exp\left(\frac{2 \, \pi \, \I}{N}\right)
  \end{equation*}
  and
  assume that
  $a$ and $k$ are positive integers satisfying $0 \leq a \leq N-1$.
  We have
  \begin{equation}
    \label{omega_sum_Bernoulli}
    \sum_{c=1}^{N-1}
    \frac{\omega_N^{~(a+1) c}}{
      \left(
        1-\omega_N^{~c}
      \right)^k
    }
    =
    \frac{(-1)^k}{(k-1)!}
    \sum_{j=1}^{k}
    \frac{S_k^{(j)}}{j} \,
    \left(
      B_{j}(1)
      -
      N^{j} \, B_{j}\left(\frac{a+1}{N}\right)
    \right)
  \end{equation}
  where $S_k^{(j)}$ is the Stirling number of the first
  kind~\eqref{generate_Stirling}.
\end{lemma}
\begin{proof}
  We follow a method in Ref.~\citen{ArakIbukKane01}.

  We define the function $P(t; k, a)$ by
  \begin{equation}
    P(t;k, a)
    =
    \sum_{c=1}^{N-1} \frac{\omega_N^{~(a+1) c}}{
      \left(
        1 - \omega_N^{~c} \, \E^t
      \right)^k}
  \end{equation}
  where $a$ and $k$ are positive integers.
  The function $P(t; k=1, a)$ is
  computed as follows;
  \begin{align}
    P(t;k=1, a)
    & =
    \sum_{n=0}^\infty \sum_{c=1}^{N-1} \omega_N^{~(a+1)c+n c} \,
    \E^{n t}
    \nonumber
    \\
    & =
    - \frac{1}{1 - \E^t}
    + N \,
    \E^{-(a+1) t}
    \, \frac{\E^{N t}}{1 - \E^{N t}}
    \nonumber
    \\
    & =
    \sum_{m=0}^{\infty}
    \frac{1}{(m+1)!} \,
    \left(
      B_{m+1}(0) -
      N^{m+1} \,
      B_{m+1}
      \left(
        1 - \frac{a+1}{N}
      \right)
    \right)
    \, t^m
    \label{solution_Pt1}
  \end{align}
  Here we have used
  \begin{equation*}
    \sum_{c=1}^{N-1} \omega_N^{~ n c}
    =
    \begin{cases}
      N-1 ,  & \text{if $N \mid n$,}
      \\[2mm]
      -1 ,  & \text{otherwise.}
    \end{cases}
  \end{equation*}
  We further introduce the function $Q(t;k, c)$ by
  \begin{equation}
    \label{define_Qt}
    Q(t; k , c)
    =
    \frac{1}{
      \left(
        1 - \omega_N^{~c} \, \E^t
      \right)^k
    }, 
  \end{equation}
  where  $k$ and $c$ are positive integers.
  By definitions we have
  \begin{equation}
    \label{Pt_and_Qt}
    P(t; k, a)
    =
    \sum_{c=1}^{N-1} \omega_N^{~ (a+1) c} \,
    Q(t ; k, c)
  \end{equation}
  The definition~\eqref{define_Qt} indicates  that
  the function $Q(t;k,c)$ satisfies a
  differential-difference
  equation,
  \begin{equation*}
    \frac{\mathrm{d}}{\mathrm{d} t}
    Q(t; k , c)
    =
    k \,
    \bigl(
      Q(t; k+1, c)
      - Q(t; k, c)
    \bigr)
  \end{equation*}
  We can check by induction that the function $Q(t; k, c)$
  can be written in terms of $Q(t; k=1, c)$ as
  \begin{equation}
    Q(t; k,c)
    =
    \frac{(-1)^{k+1}}{(k-1)!} 
    \sum_{m=0}^{k-1}
    (-1)^m \, S_k^{(m+1)} \,
    \frac{\mathrm{d}^m}{\mathrm{d} t^m}
    Q(t;1,c)
  \end{equation}
  From~\eqref{Pt_and_Qt}
  we find  that $P(t;k,a)$ is solved as
  \begin{equation}
    P(t;k,a)
    =
    \frac{(-1)^{k+1}}{(k-1)!}
    \sum_{m=0}^{k-1}
    (-1)^m \, S_k^{(m+1)} \,
    \frac{\mathrm{d}^m}{\mathrm{d} t^m}
    P(t;1,a)
  \end{equation}
  Substituting~\eqref{solution_Pt1} for the above solution,
  we complete the proof.
\end{proof}

Using the arithmetic identity~\eqref{omega_sum_Bernoulli}, we can
rewrite the WRT
invariant~\eqref{result_Rozansky} in terms of the Bernoulli
polynomials.

\begin{prop}
  \label{prop:invariant_Bernoulli}
  The WRT invariant for the $M$-exceptional fibered Seifert integral
  homology
  sphere $\mathcal{M}=\Sigma(\vv{p})$, which was computed as
  in~\eqref{result_Rozansky},
  is written  in terms of the Bernoulli polynomials as
%
    \begin{multline}
      \label{invariant_Bernoulli_large}
      \E^{\frac{2 \pi \I}{N}
        \left( \frac{\phi(\vv{p})}{4} - \frac{1}{2} \right)
      } \,
      \left(
        \E^{\frac{2 \pi \I}{N}} - 1
      \right) \,
      \tau_N
      (\mathcal{M})
      =
      -\frac{1}{2} \, \frac{1}{(M-3)!} \sum_{j=0}^{M-3} (-N)^j \,
      \frac{S_{M-2}^{(j+1)}}{j+1}
      \\
      \times
      \sum_{n=0}^{2 P N -1} 
      \chi_{2 P}^{\vv{E}}(n) \,
      \E^{\frac{1}{2 P N} \left( n+ P (M-3) \right)^2 \pi
        \I} \,
      B_{j+1}
      \left(
        \frac{1}{N} \, \left\lfloor \frac{n}{2 \, P}\right\rfloor
      \right)
      \\
      +
      \frac{(-1)^{M}}{(M-4)!}
      \sum_{a=1}^\infty \sum_{b=0}^{a-1}
      \sideset{}{^\prime}\sum_{\vv{\eta}(a)}
    \sum_{n=0}^{N-1} \sum_{j=1}^{M-3}
    \frac{S_{M-3}^{(j)}}{j}
    \,
    \E^{\pi \I \frac{2 P}{N}
      \left(
        n - b - \frac{M-4}{2} + \frac{1}{2}
        \sum_{i=1}^M \frac{\eta_i}{p_i}
      \right)^2
    } \,
    \\
    \times
    \left(
      N^{j-1} \,
      B_j\left(
        \frac{n+1}{N}
      \right)
      - \frac{1}{N} \, B_{j}(1)
    \right)
  \end{multline}
  where for our brevity we have used $M$-tuple
    \begin{equation}
      \vv{E}=
      (\underbrace{1,1,\dots, 1}_{M})
    \end{equation}
    and  $\displaystyle\sideset{}{^\prime}\sum_{\vv{\eta}(a)}$ is a signed sum
    \begin{equation*}
      \sideset{}{^\prime}\sum_{\vv{\eta}(a)}
      \cdots
      =
      \sum_{
        \substack{
          \vv{\eta}\in \{ \pm 1\}^{ M}
          \text{s.t.}
          \\
          2 a -1 < \sum_{j=1}^M \frac{\eta_j}{p_j} < 2 a +1}}
      \left[
        \prod_{j=1}^M \eta_j
      \right]
      \cdots
    \end{equation*}
\end{prop}

We remark that
the second term including a sum of $a$ in~\eqref{invariant_Bernoulli_large}
vanishes when
\begin{equation*}
  \sum_{j=1}^M \frac{1}{p_j} <1
\end{equation*}
Even when $\sum_{j=1}^M \frac{1}{p_j} >1$, the second term is a
finite sum.
It is well known that the sum of inverse of prime numbers,
$\sum_{\text{$p$: prime}} \frac{1}{p}$ diverges, although
the sum up to the $10,000$-th prime numbers is still $2.709258 \cdots$.

\begin{proof}[Proof of Prof.~\ref{prop:invariant_Bernoulli}]
  We first study a case of $\sum_j \frac{1}{p_j}<1$.
  In this case, we have
  \begin{equation}
    \label{product_and_chi}
    - z^P \,
    \prod_{j=1}^M
    \left(
      z^{\frac{P}{p_j}} - z^{- \frac{P}{p_j}}
    \right)
    =
    \sum_{m=0}^{2 P -1}
    \chi_{2 P}^{\vv{E}}(m) \, z^m
  \end{equation}
  where the periodic function $\chi_{2 P}^{\vv{E}}(m)$ is defined
  in~\eqref{define_chi}.
  Using this identity in~\eqref{result_Rozansky},
  we  have
  \begin{align*}
    &
    \E^{\frac{2 \pi \I}{N}
      \left( \frac{\phi(\vv{p})}{4} - \frac{1}{2} \right)
    } \,
    \left(
      \E^{\frac{2 \pi \I}{N}} - 1
    \right) \,
    \tau_N
    (\mathcal{M})
    \\
    & =
    - \frac{\E^{\frac{\pi \I}{4}}}{2
      \,
      \sqrt{2 \, P \, N}
    }
    \sum_{m=0}^{2 P -1}
    \sum_{k=1}^{N-1}
    \chi_{2 P}^{\vv{E}}(m) \,
    \,
    \frac{
      \E^{\frac{\pi \I}{2 P N}
        \left(
          m + (M-3) P
        \right)^2
      }}{
      \left(
        \E^{\frac{2 \pi \I}{N} k} - 1
      \right)^{M-2}
    }
    \sum_{j =0}^{2 P -1}
    \E^{-
      \pi \I \frac{N}{2 P}
      \left(
        j+
        \frac{ k - m - (M-3) P}{N}
      \right)^2
    }
    \\
    & =
    \frac{-1}{2 \, N}
    \sum_{m=0}^{2 P -1} \sum_{n=0}^{N-1}
    \chi_{2 P}^{\vv{E}}(m) \,
    \E^{\frac{\pi \I}{2 P N}
      \left(
        2 P n - 
        m - (M-3) P
      \right)^2
    }
    \sum_{k=1}^{N-1}
    \frac{
      \E^{\frac{2 \pi \I}{N} k n}
    }{
      \left(
        \E^{\frac{2 \pi \I k}{N}} - 1
      \right)^{M-2}
    }
    \\
    & =
    \frac{-1}{2} \,
    \frac{1}{(M-3)!}
    \sum_{m=0}^{2 P -1} \sum_{n=0}^{N-1}
    \chi_{2 P}^{\vv{E}}(m) \,
    \E^{
      \frac{\pi \I}{2 P N}
      \left( 2 P n + m + (M-3) P \right)^2
    }
    \\
    & \qquad \qquad \qquad
    \times
    \sum_{j=0}^{M-3}
    \frac{S_{M-2}^{(j+1)}}{j+1}
    \left(
      \frac{1}{N} \, B_{j+1}(1)
      -
      N^j \,
      B_{j+1}
      \left(
        1 - \frac{n}{N}
      \right)
    \right)
  \end{align*}
  Here in the first equality we have used~\eqref{product_and_chi},
  and decomposed  a sum of $n$ by setting $n = N \, j + k$.
  We have then applied the Gauss sum reciprocity
  formula~\eqref{Gauss_reciprocity} in the second equality,
  and then used our formula~\eqref{omega_sum_Bernoulli} in the last
  equality.
  As we  have
  \begin{equation}
    \label{sum_limit_zero}
    \sum_{n=0}^{2 P N-1}
    \chi_{2 P}^{\vv{\ell}}(n) \,
    \E^{\frac{\pi \I}{2 P N}
      n^2
    } = 0
  \end{equation}
  due to~\eqref{chi_mean_zero},
  we obtain the first term of~\eqref{invariant_Bernoulli_large}.

  For other cases $\sum_j \frac{1}{p_j} > 1$,
  the generating function~\eqref{product_and_chi} of the periodic
  function $\chi_{2 P}^{\vv{E}}(n)$ is replaced with
  \begin{multline}
    \label{product_and_chi_others}
    - z^P \,
    \prod_{j=1}^M
    \left(
      z^{\frac{P}{p_j}} - z^{- \frac{P}{p_j}}
    \right)
    +
    \sum_{a=1}^\infty
    \sideset{}{^\prime}\sum_{\vv{\eta}(a)}
    z^P \,
    \left(
      z^{a P} - z^{- aP}
    \right)
    \\
    \times
    \left(
      z^{P \left( \sum_j \frac{\eta_j}{p_j} - a
        \right)}
      +
      (-1)^{M+1} \,
      z^{ - P \left( \sum_j \frac{\eta_j}{p_j} - a
        \right)}
    \right)
    = \sum_{n=0}^{2 P -1}
    \chi_{2 P}^{\vv{E}}(n) \, z^n
  \end{multline}
  Thus comparing with a case of $\sum_j\frac{1}{p_j}<1$, we need
  additional term $\tau_{\text{add}}$ defined by
  \begin{multline}
    \tau_{\text{add}} =
    \sum_{a=1}^\infty \sum_{b=0}^{a-1}
    \sideset{}{^\prime}\sum_{\vv{\eta}(a)}
    \sum_{
      \substack{
        n=0
        \\
        N \nmid n
      }}^{2 P N -1}
    \frac{\E^{\frac{\pi \I}{4}}}{2 \sqrt{2 \, P \, N}}
    \frac{
      \E^{-\frac{n^2}{2 P N} \pi \I}
    }{
      \left(
        \E^{\frac{n}{N} \pi \I} -
        \E^{-\frac{n}{N} \pi \I}
      \right)^{M-3}
    }
    \\
    \times
    \left(
      \E^{ \pi \I \frac{n}{N}
        \left(
          \sum_j \frac{\eta_j}{p_j}-2 b -1
        \right)
      }
      +(-1)^{M+1} \,
      \E^{- \pi \I \frac{n}{N}
        \left(
          \sum_j \frac{\eta_j}{p_j}-2 b -1
        \right)
      }
    \right)
  \end{multline}
  We decompose a sum of $n$ by setting $n=N \, j+k$, and apply the
  Gauss sum reciprocity formula~\eqref{Gauss_reciprocity}.
  After some computations, we obtain
  \begin{align*}
    \tau_{\text{add}}
    &=
    \frac{1}{2 \, N}
    \sum_{a=1}^\infty \sum_{b=0}^{a-1}
    \sideset{}{^\prime}\sum_{\vv{\eta}(a)}
    \sum_{k=1}^{N-1} \sum_{n=0}^{N-1}
    \frac{
      \E^{\pi \I \frac{2 P}{N}
        \left(
          n - \frac{M-2}{2} + \frac{1}{2} \sum_j \frac{\eta_j}{p_j}
        \right)^2
      }
    }{
      \left(
        \E^{\frac{k}{N} \pi \I} -
        \E^{-\frac{k}{N} \pi \I}
      \right)^{M-3}
    }
     \\
     & \qquad \quad
     \times
    \left(
      \E^{\pi \I \frac{k}{N} (M-3-2b-2n)}
      +(-1)^{M+1}
      \E^{\pi \I \frac{k}{N} (2 b + 2 n - M+3)}
    \right)
    \\
    & =
    \frac{1}{N}
    \sum_{a=1}^\infty \sum_{b=0}^{a-1}
    \sideset{}{^\prime}\sum_{\vv{\eta}(a)}
    \sum_{n=0}^{N-1}
    \E^{\pi \I \frac{2P}{N}
      \left(
        n - \frac{M-2}{2} + \frac{1}{2} \sum_j \frac{\eta_j}{p_j}
      \right)^2
    }
    \sum_{k=1}^{N-1}
    \frac{
      \E^{2 \pi \I \frac{k}{N} (b+n)}
    }{
      \left(
        1 - \E^{2 \pi \I \frac{k}{N}}
      \right)^{M-3}
    }
  \end{align*}
  where we have used the fact that the sum~\eqref{omega_sum_Bernoulli}
  is real.
  Substituting~\eqref{omega_sum_Bernoulli} for the above, we
  find that $\tau_{\text{add}}$ gives the second term of
  r.h.s. of~\eqref{invariant_Bernoulli_large}, and thus we complete
  the proof.
\end{proof}

We now aim  to relate this expression with limiting values of
the Eichler integrals~\eqref{limiting_value_Eichler} and
differentials~\eqref{limiting_frac_differential}  thereof.
For our convention, we introduce an analogue of the Bernoulli polynomial
defined by
\begin{equation}
  \label{define_poly}
  f_m^M(x)
  =
  \sum_{k=m}^{M}
  \frac{1}{k} \,
  S_M^{(k)} \,
  \begin{pmatrix}
    k \\
    m
  \end{pmatrix}
  \,
  \left(
    x+ \frac{M}{2}
  \right)^{k-m}
\end{equation}
where $M, m \in\mathbb{Z}$ satisfying $M \geq m>0$.
For a case of $m=0$, we set
\begin{equation}
  f_0^M(x)
  =
  \sum_{k=1}^M \frac{S_M^{(k)}}{k} \,
  \left(
    x+ \frac{M}{2} 
  \right)^k
\end{equation}
Some of  explicit forms of the polynomials $f_m^M(x)$ are given below;
\begin{align*}
  f_M^M(x)
  & = \frac{1}{M} \\[2mm]
  f_{M-1}^M(x) & = x
  \\[2mm]
  f_{M-2}^M(x) & =
  \frac{1}{M} \,
  \begin{pmatrix}
    M \\ 2
  \end{pmatrix} \,
  \left(
    x^2 - \frac{M}{12}
  \right)
  \\[2mm]
  f_{M-3}^M(x) & =
  \frac{1}{M} \,
  \begin{pmatrix}
    M \\ 3
  \end{pmatrix} \,
  \left(
    x^3  - \frac{M}{4} \, x
  \right)
  \\[2mm]
  f_{M-4}^M(x) & =
  \frac{1}{M} \,
  \begin{pmatrix}
    M \\ 4
  \end{pmatrix} \,
  \left(
    x^4 - \frac{M}{2} \, x^2 +
    \frac{1}{240} \, M \, (5 \, M+2)
  \right)
\end{align*}

\begin{lemma}
  \label{lemma:f_even_odd}
  Let the polynomial $f_m^M(x)$ be defined by~\eqref{define_poly}.
  Then the polynomial  $f_{M-k \neq 0}^M(x)$ is even
  (resp. odd)
  when
  $k$ is
  even (resp. odd).
\end{lemma}

\begin{proof}
  We introduce  the generating function of the polynomials $f_M^m(x)$ by
  \begin{equation}
    F_M(x,y) = \sum_{m=0}^M m \, f_m^M (x) \, y^{m-1}
  \end{equation}
  Recalling the generating function~\eqref{generate_Stirling}
  of the Stirling number of the
  first kind, we get
  \begin{align}
    \label{generate_f_poly}
    F_M(x,y)
    & =
    \prod_{j=1}^{M-1} \left( y + x+\frac{M}{2} - j \right)
  \end{align}
  which shows that $F_M(x,y)$ is a polynomial of $x+y$.
  Furthermore  $F_M(x,y)$ becomes an odd (resp. even)
  polynomial of $x+y$ when $M$ is even (resp. odd).
  Then  we can conclude that
  the polynomial $f_{M-k}^M(x)$ is even (resp. odd)
  if $k$ is even (resp. odd).
\end{proof}

By use of the generating function~\eqref{generate_f_poly},
we obtain the following differential equation and recursion relation of $f_m^M(x)$;
\begin{gather}
  \frac{
    \mathrm{d}}{\mathrm{d} x}
  f_m^M(x)
  =
  (m+1) \, f_{m+1}^M(x)
  \\[2mm]
  \label{recursion_f_poly}
  f_j^{M+1}
  \left( x- \frac{1}{2} \right)
  =
  \left(
    x- \frac{M}{2}
  \right) \,
  f_j^M
  \left(
    x \right)
  +\frac{j-1}{j} \,
  f_{j-1}^M  \left(x \right)
\end{gather}
We can rewrite the WRT invariant in
Prop.~\ref{prop:invariant_Bernoulli} in terms of these polynomials as
follows.

\begin{prop}
  The WRT invariant $\tau_N(\mathcal{M})$
  for the $M$-exceptional fibered Seifert homology sphere
  $\mathcal{M}=\Sigma(\vv{p})$
  is written as
%
    \begin{multline}
      \label{invariant_Bernoulli_2_large}
      \E^{\frac{2 \pi \I}{N}
        \left( \frac{\phi(\vv{p})}{4} - \frac{1}{2} \right)
      } \,
      \left(
        \E^{\frac{2 \pi \I}{N}} - 1
      \right) \,
      \tau_N
      (\mathcal{M})
      \\
      =
      -\frac{1}{2} \, \frac{1}{(M-3)!}
      \sum_{k=1}^{M-2} (-N)^{k-1} \,
      \sum_{n=0}^{2 P N -1}
      \chi_{2 P}^{\vv{E}}(n) \,
      f_k^{M-2}\left(
        \frac{n}{2 \, P} -
        \left\lfloor
          \frac{n}{2 \, P}
        \right\rfloor
        - \frac{1}{2}
      \right)
      \\
      \times
      \E^{\frac{1}{2 P N}
        \left(
          n+ P (M-3)
        \right)^2 \pi \I
      } \,
      B_k
      \left(
        \frac{n+ P \, (M-3)}{2 \, P \, N}
      \right)
      \\
      +
      \frac{(-1)^{M}}{(M-4)!}
      \sum_{a=1}^\infty \sum_{b=0}^{a-1}
      \sideset{}{^\prime}\sum_{\vv{\eta}(a)}
      \sum_{n=0}^{N-1}
      \E^{\pi \I \frac{2 P}{N}
        \left( 
          n - b - \frac{M-4}{2} + \frac{1}{2} \sum_j
          \frac{\eta_j}{p_j}
        \right)^2
      }
      \\
      \times
      \sum_{k=0}^{M-3}
      \Biggl(
        N^{k-1} \, 
        f_k^{M-3}\left(
          b + \frac{1}{2} - \frac{1}{2} \sum_j \frac{\eta_j}{p_j}
        \right) \,
        B_k\left(
          \frac{
            2 \, n - 2 \, b - M+4 + \sum_j \frac{\eta_j}{p_j}
          }{
            2 \, N
          }
        \right)
        \\
        -
        \frac{1}{N}
        f_k^{M-3}\left(
          \frac{5-M}{2}
        \right) \,
        B_k(0)
      \Biggr)
    \end{multline}
\end{prop}
\begin{proof}
  To prove  for a  case of $\sum_j \frac{1}{p_j} <1$,
  in which the second term in~\eqref{invariant_Bernoulli_2_large}
  vanishes,
  we only need to apply~\eqref{Bernoulli_x_y} to
  $B_{j+1}
  \left(
    \frac{1}{N}
    \left\lfloor
      \frac{n}{2 P}
    \right\rfloor
  \right)$
  in the first term of~\eqref{invariant_Bernoulli_large}.
  As we have an identity
  \begin{equation}
    \label{f_poly_Bernoulli}
    \sum_{j=1}^M \frac{S_M^{(j)}}{j} \, B_j(x+y) \, z^j
    =
    \sum_{k=0}^M f_k^M \left(y \, z - \frac{M}{2}\right) \, B_k(x) \,
    z^k
  \end{equation}
  we obtain the required expression.

  In the case of $\sum_j \frac{1}{p_j} > 1$,
  we 
  need to evaluate the
  second term in~\eqref{invariant_Bernoulli_large}, which we have set
  $\tau_{\text{add}}$ in the proof of
  Prop.~\ref{prop:invariant_Bernoulli}.
  This term $\tau_{\text{add}}$ can be transformed into the above
  expression~\eqref{invariant_Bernoulli_2_large}
  when we apply~\eqref{f_poly_Bernoulli}.
\end{proof}

\begin{theorem}
  The exact asymptotic expansion of the WRT invariant
  $\tau_N(\mathcal{M})$
  for the Seifert
  homology sphere
  $\mathcal{M}=\Sigma(\vv{p})$ with $M$-singular fibers
  in $N\to\infty$ is given as follows;
    \begin{multline}
      \label{exact_asymptotics_large_1}
      \E^{\frac{2 \pi \I}{N}
        \left( \frac{\phi(\vv{p})}{4} - \frac{1}{2} \right)
      } \,
      \left(
        \E^{\frac{2 \pi \I}{N}} - 1
      \right) \,
      \tau_N
      (\mathcal{M})
      \\
      \simeq
      \frac{1}{2}  \, \frac{1}{(M-3)!}  \,
      \frac{1
      }{\I^{\frac{3}{2} - \eveodd(M)}}
      \sum_{j=0}^{
        \left\lfloor \frac{M-3}{2} \right\rfloor
      }
      N^{j+
        \left\lfloor
          \frac{M}{2}
        \right\rfloor
        - \frac{1}{2}} 
      \,
      \left(
        \frac{\I}{2 \, P \,  \pi }
      \right)^{
        \lfloor \frac{M-3}{2} \rfloor - j
      }
      \,
      K_{\left\lfloor \frac{M-3}{2} \right\rfloor,
        \eveodd(M)}^{(j)}
      \,
      \\
      \times
      \frac{1}{
        2 \, j+2 - \eveodd(M)
      } 
      \sum_{\vv{\ell^\prime}}
      \mathbf{S}^{\sigma_1^{M-1}(\vv{E})}_{\vv{\ell^\prime}} \,
      \left[
        \sum_{n=1}^{2 P}
        \chi_{2 P}^{\vv{\ell^\prime}}(n) \,
        B_{2j+2 - \eveodd(M)}
        \left(
          \frac{n}{2 P}
        \right)
      \right] \,
      \E^{- N \frac{P}{2}
        \left(
          1+ \sum_j \frac{\ell_j^\prime}{p_j}
        \right)^2 \pi \I
      }
      \\
      +
      \frac{1}{(M-3)!} \,
      \frac{(-1)^M}{\I^{\frac{3}{2}}}
      \sum_{m=1}^{
        \left\lfloor \frac{M-3}{2} \right\rfloor}
      \sum_{a=0}^P
      \chi_{2 P}^{\sigma_1^{M-1}(\vv{E})}(a) \,
      f_{2 m}^{M-2}
      \left(
        \frac{a}{2 \, P} - \frac{\eveodd(M)}{2}
      \right)
      \\
      \times
      \sum_{j=1}^{m}
      N^{m+j-\frac{1}{2}} \,
      \left(
        \frac{
          \I
        }{
          2 \, P  \, \pi  
        }
      \right)^{m-j}
      \,  K_{m-1,0}^{(j-1)} \,
      \frac{m}{j} \,
      \sum_{c=0}^P
      \mathbf{N}^a_c \,
      \frac{2 - \delta_{c,0} - \delta_{c,P}}{2} \,
      B_{2j}
      \left(
        \frac{c}{2 \, P}
      \right) \,
      \E^{- N \frac{c^2}{2 P} \pi \I}
      \\
      +
      \frac{1}{(M-3)!} \,
      \frac{(-1)^{M-1}}{
        \I^{\frac{1}{2}}
      }
      \sum_{m=0}^{
        \left\lfloor \frac{M-4}{2} \right\rfloor
      }
      \sum_{a=1}^P
      \chi_{2 P}^{\sigma_1^{M-1}(\vv{E})}(a) \,
      f_{2m+1}^{M-2}
      \left(
        \frac{a}{2 \, P} - \frac{\eveodd(M)}{2}
      \right)
      \\
      \times
      \sum_{j=0}^m
      N^{m +j + \frac{1}{2} }    \,
      \left(
        \frac{\I}{2  \, P  \, \pi }
      \right)^{m-j}
      \, K_{m,1}^{(j)} \,
      \frac{2 \, m+1}{2 \, j +1}
      \sum_{c=1}^{P-1}
      \mathbf{M}^{a}_c
      \, B_{2j+1}
      \left(
        \frac{c}{2 \, P}
      \right) \,
      \E^{-N \frac{c^2}{2 P} \pi \I}
      \\
      +
      \frac{(-1)^M}{(M-4)!} \sum_{a=1}^\infty
      \sum_{b=0}^{a-1}
      \sideset{}{^\prime}\sum_{\vv{\eta}(a)}
      \Biggl\{
        \sum_{m=1}^{\lfloor \frac{M-3}{2} \rfloor}
        f_{2m}^{M-3}\left(b+\frac{1}{2}-\frac{1}{2}\sum_i
          \frac{\eta_i}{p_i}
        \right)
        \,
        \frac{-1}{\I^{\frac{3}{2}}} \sum_{j=1}^{m} N^{m+j-\frac{1}{2}}
        \left(
          \frac{\I}{2 \, P \, \pi }
        \right)^{m-j}
        \\
        \times
        K_{m-1,0}^{(j-1)}
        \frac{m}{j}
        \sum_{c=0}^P
        \mathbf{N}^{\sigma^{M-1}(P\sum_i\frac{\eta_i}{p_i}-P(2a-1))}_{c}
        \,
        \frac{2 - \delta_{c,0} - \delta_{c,P}}{2} \,
        B_{2j}\left(\frac{c}{2 P}\right) \,
        \E^{-N  \frac{c^2}{2 P} \pi \I}
        \\
        +(-1)^{M-1}
        \sum_{m=0}^{\lfloor \frac{M-4}{2} \rfloor}
        f_{2m+1}^{M-3}\left(b+\frac{1}{2}-\frac{1}{2}\sum_i
          \frac{\eta_i}{p_i}
        \right)
       ,
        \frac{-1}{\I^{\frac{1}{2}}} \sum_{j=0}^{m} N^{m+j+\frac{1}{2}}
        \left(
          \frac{\I}{2 \, P \, \pi }
        \right)^{m-j}
        \\
        \times
        K_{m,1}^{(j)}
        \frac{2 \, m+1}{2 \, j+1}
        \sum_{c=1}^{P-1}
        \mathbf{M}^{\sigma^{M-1}(P\sum_i\frac{\eta_i}{p_i}-P(2a-1))}_{c}
        B_{2j+1}\left(\frac{c}{2 P}\right) \,
        \E^{-N  \frac{c^2}{2 P} \pi \I}
      \Biggr\}
      \\
      +
      \sum_{k=0}^\infty
      \frac{T_{\vv{p}}(k)}{k!} \,
      \left(
        \frac{
          \pi \, \I
        }{
          2 \, P \, N
        }
      \right)^k
    \end{multline}
    where we have used  an involution $\sigma$ on $x\in \mathbb{Z}_{2 P}$
    defined by
    \begin{equation}
      \sigma(x) = P - x \mod 2 \, P
    \end{equation}
    The coefficients $T_{\vv{p}}(k)$ in a \emph{tail} part are defined by
    \begin{multline}
      \label{define_T_large_1}
      T_{\vv{p}}(k)
      =
      \frac{(-1)^{M+1}}{2} \, \frac{(2 \, P)^{2 k}}{(M-3)!}
      \sum_{n=1}^{2 P}
      \chi_{2 P}^{\sigma_1^{M-1}(\vv{E})}(n)
      \sum_{j=1}^{M-2}
      (-1)^j\,
      \frac{j}{2 \, k+j} \,
      B_{2k+j}
      \left(
        \frac{n}{2 \, P}
      \right)
      \\
      \times
      f_j^{M-2}
      \left(
        \frac{n+(M-1) \, P}{2 \, P}
        -\left\lfloor
          \frac{n+(M-1) \, P}{2 \, P}
        \right\rfloor
        -\frac{1}{2}
      \right)
      \\
      -      
      \frac{(2 \,P)^{2k}}{(M-4)!} \sum_{a=1}^\infty \sum_{b=0}^{a-1}
      \sideset{}{^\prime}\sum_{\vv{\eta}(a)}
      \sum_{j=1}^{M-3}
      (-1)^{(M+1)(j+1)} \, \frac{j}{2 \, k +j}
      \\
      \times
      f_j^{M-3}\left(
        b +\frac{1}{2} - \frac{1}{2} \sum_i \frac{\eta_i}{p_i}
      \right) \,
      B_{2k+j}\left(
        \frac{1}{2 \, P} \,
        \sigma^{M-1}\left(
          P\sum_i\frac{\eta_i}{p_i} - (2 \, a -1) \, P
        \right)
      \right)
    \end{multline}
\end{theorem}
\begin{proof}
  We first study the case $\sum_j \frac{1}{p_j}<1$,
  in which we only have the first term
  in~\eqref{invariant_Bernoulli_2_large}.
  We shift the parameter $n$ by $P \, (3-M)$, and
  we have
  \begin{equation*}
    \chi_{2P}^{\vv{E}}(n-P \, (M-3)) =
    (-1)^{1 - \eveodd(M)} \, \chi_{2P}^{\sigma_1^{M-1}(\vv{E})}(n)
  \end{equation*}
  As 
  we see that
  $
  \chi_{2 P}^{\vv{\ell}}(n) \,
  f_{k}^{M-2}
  \left(
    \frac{n+(M-1) P}{2 \,P}
    -
    \left\lfloor
      \frac{n+(M-1) P}{2 \,P}
    \right\rfloor
    -
    \frac{1}{2}
  \right)
  $ is an even (resp. odd) periodic function when $k$ is even (resp. odd),
  we can write by use of  $f_M^M(x)=\frac{1}{M}$ that
  \begin{multline*}
    \E^{\frac{2 \pi \I}{N}
      \left( \frac{\phi(\vv{p})}{4} - \frac{1}{2} \right)
    } \,
    \left(
      \E^{\frac{2 \pi \I}{N}} - 1
    \right) \,
    \tau_N
    (\mathcal{M})
    \\
    =
    -\frac{1}{2} \, \frac{N^{M-3}}{(M-2)!} \,
    \sum_{n=0}^{2 P N -1}
    \chi_{2 P}^{\sigma_1^{M-1}(\vv{E})}(n) \,
    \E^{\frac{n^2}{2 P N} \pi \I} \,
    B_{M-2}\left(\frac{n}{2 \, P \, N}\right)
    \\
    - \frac{(-1)^M}{2} \, \frac{1}{(M-3)!} \,
    \sum_{c=1}^{\lfloor \frac{M-3}{2} \rfloor}
    N^{2 c-1}    \,
    \left[
      \sum_{n=0}^{2 P N -1}
      \chi_{2 P}^{\sigma_1^{M-1}(\vv{E})}(n) \,
      f_{2c}^{M-2}\left(\frac{n}{2 \, P} - \frac{\eveodd(M)}{2}\right)
      \,
      \E^{\frac{n^2}{2 P N} \pi \I} \,
      B_{2c} \left(\frac{n}{2 \, P \, N}\right)
    \right]
    \\
    + \frac{(-1)^M}{2} \, \frac{1}{(M-3)!} \,
    \sum_{c=0}^{\lfloor \frac{M-4}{2} \rfloor}
    N^{2 c} \,
    \left[
      \sum_{n=0}^{2 P N -1}
      \chi_{2 P}^{\sigma_1^{M-1}(\vv{E})}(n) \,
      f_{2c+1}^{M-2}\left(\frac{n}{2 \, P} - \frac{\eveodd(M)}{2}\right)
      \,
      \E^{\frac{n^2}{2 P N} \pi \I} \,
      B_{2c+1} \left(\frac{n}{2 \, P \, N}\right)
    \right]
  \end{multline*}
  Generating function~\eqref{product_and_chi} proves
  \begin{equation*}
    \sum_{n=0}^{2 P-1} \chi_{2 P}^{\vv{E}}(n) \,
    g(n) =0
  \end{equation*}
  where  $g(n)$ is an arbitrary polynomial of $n$ of order at most $M-1$.
  So we find that
  $      \chi_{2 P}^{\vv{E}}(n) \,
  f_k^{M-2}\left(
    \frac{n}{2 \, P} -
    \left\lfloor
      \frac{n}{2 \, P}
    \right\rfloor
    - \frac{1}{2}
  \right)
  $ is a periodic function of $n$ with mean value zero, and that
  the
  expression~\eqref{invariant_Bernoulli_2_large} can be identified with
  a limiting value of
  the Eichler integrals and their derivatives studied in the previous
  section.
  This proves that the WRT invariant is a limiting value of the
  holomorphic function of $q$ in $q\to \E^{2 \pi \I/N}$~\cite{KHabi02a}.
  Substituting both~\eqref{asymptotic_psi_Bernoulli}
  and~\eqref{asymptotic_theta_Bernoulli}
  for above expression,
  we get~\eqref{exact_asymptotics_large_1}.
  
  For the second term $\tau_{\text{add}}$
  in~\eqref{invariant_Bernoulli_2_large},
  recalling the periodic functions
  $\psi_{2P}^{(a)}(n)$ and $\theta_{2P}^{(a)}(n)$
  we can rewrite it into
  \begin{multline*}
    \tau_{\text{add}}
    =
    \frac{1}{2} \frac{(-1)^M}{(M-4)!} \,
    \sum_{a=1}^\infty \sum_{b=0}^{a-1}
    \sideset{}{^\prime}\sum_{\vv{\eta}(a)}
    \Biggl\{
      \sum_{m=0}^{\lfloor \frac{M-3}{2} \rfloor}
      N^{2m-1} \,
      f_{2m}^{M-3}\left(
        b+\frac{1}{2} - \frac{1}{2} \sum_i \frac{\eta_i}{p_i}
      \right)
      \\
      \times
      \sum_{n=0}^{2PN-1}
      \theta_{2P}^{(\sigma^{M-1}(P\sum_i\frac{\eta_i}{p_i} -(2a-1)P))}(n) \,
      \E^{\pi \I \frac{n^2}{2 P N}} \,
      B_{2m}\left(\frac{n}{2 \, P \, N}\right)
      \\
      +(-1)^{M-1} 
      \sum_{m=0}^{\lfloor \frac{M-4}{2} \rfloor}
      N^{2m} \,
      f_{2m+1}^{M-3}\left(
        b+\frac{1}{2} - \frac{1}{2} \sum_i \frac{\eta_i}{p_i}
      \right)
      \\
      \times
      \sum_{n=0}^{2PN-1}
      \psi_{2P}^{(\sigma^{M-1}(P\sum_i\frac{\eta_i}{p_i} -(2a-1)P))}(n) \,
      \E^{\pi \I \frac{n^2}{2 P N}} \,
      B_{2m+1}\left(\frac{n}{2 \, P \, N}\right)
    \Biggr\}
  \end{multline*}
  With this term, the mean value zero condition is satisfied even in
  this case, and we can apply
  the result of Coro.~\ref{coro:asymptotic_Bernoulli}
  to obtain the exact asymptotic expansion as required.
\end{proof}

We have thus obtained that the WRT invariant $Z_{N-2}(\mathcal{M})$
for the  $M$-exceptional fibered Seifert homology sphere
$\mathcal{M}=\Sigma(\vv{p})$
is written in a limit $N\to \infty$
as a sum of exponentially divergent terms and a \emph{tail} part;
\begin{equation}
  \label{Z_sum_decompose}
  Z_{N-2}(\mathcal{M})
  \simeq
  \sum_{k=0}^{M-3}
  N^{M-3-k} \,
  Z_{N-2}^{(k)}(\mathcal{M})
  +
  \text{a tail part}
\end{equation}
Here a \emph{tail part} means an infinite power series of $1/N$,
and it corresponds to a contribution from the trivial connection.
Among the divergent terms in~\eqref{Z_sum_decompose},
the dominating term in the limit $N\to\infty$
is $Z_{N-2}^{(0)}(\mathcal{M})$, which is  read as follows.

\begin{coro}
  In a limit $N\to \infty$, 
  the asymptotics of the WRT invariant $\tau_N(\mathcal{M})$ for the
  $M$-exceptional fibered Seifert homology sphere $\mathcal{M}=\Sigma(\vv{p})$
  is dominated by $N^{M-3} \cdot Z_{N-2}^{(0)}(\mathcal{M})$, namely
  \begin{multline}
    \label{WRT_dominating}
    Z_{N-2}(\mathcal{M})
    \sim
    N^{M-3} \cdot
    \E^{- \frac{\phi(\vv{p})}{2 N} \pi \I}
    \,
    \frac{
      \I^{\eveodd(M) -1} \,
      \E^{-\frac{3}{4} \pi \I}
    }{
      2 \, \sqrt{2} \, (M-2)!
    }
    \\
    \times
    \sum_{\vv{\ell}}
    \mathbf{S}^{\sigma_1^{M-1}(\vv{E})}_{\vv{\ell}} 
    \,
    C_{\vv{p}} \left( \vv{\ell} \right)
    \,
    \E^{-  \frac{P}{2} 
      \left(
        1+ \sum_{j=1}^M \frac{\ell_j}{p_j}
      \right)^2 \pi \I N
    } .
  \end{multline}
  \emph{i.e.},
  \begin{multline}
    Z_{N-2}(\mathcal{M})
    \sim
    N^{M-3} \, \frac{2^{M-2}}{(M-2)! \, \sqrt{P}} \,
    \E^{-\frac{\phi(\vv{p})}{2 N} \pi \I} \,
    \E^{-\frac{2 M+1}{4} \pi \I}
    \sum_{\vv{\ell}}
    C_{\vv{p}} \left( \vv{\ell} \right)
    \, \E^{-\frac{P}{2} \left( 1 + \sum_{j=1}^M
        \frac{\ell_j}{p_j} \right)^2 \pi \I N}
    \\
    \times
    (-1)^{M P \left( 1+\sum_j \frac{\ell_j}{p_j} \right)
      +P \sum_j \frac{1}{p_j} + P \sum_j \sum_{k \neq j}
      \frac{\ell_k}{p_j p_k}
    } \,
    \left[
      \prod_{j=1}^M
      \sin \left(
        P \, \frac{\ell_j}{p_j^{~2}} \, \pi
      \right)
    \right]
  \end{multline}
  where the sum of $M$-tuples
  $\vv{\ell}$ runs over $D$-dimensional
  space~\eqref{define_Dp},
  and
  the function $C_{\vv{p}} \left( \vv{\ell} \right)$ is defined by
  \begin{equation}
    \label{define_Cpl}
    C_{\vv{p}} \left( \vv{\ell} \right)
    =
    \sum_{n=1}^{2 \, P}
    \chi_{2 P}^{\vv{\ell}}(n) \,
    B_{M-2}
    \left(
      \frac{n}{2 \, P}
    \right)
  \end{equation}
\end{coro}

We note that
the invariance~\eqref{sigma2_chi} of the periodic function
$\chi_{2P}^{\vv{\ell}}(n)$ indicates that
\begin{equation}
  \label{C_polynomial_invariant}
  C_{\vv{p}} \left( \sigma_i \sigma_j (\vv{\ell}) \right)
  =
  C_{\vv{p}} \left( \vv{\ell} \right)
\end{equation}

By construction, the asymptotics~\eqref{WRT_dominating} should coincide
with~\eqref{asymptotics_Lawrence}  which follows from residue part
of~\eqref{tau_res_int}.
We do not have a direct proof,  and we have checked the equivalence
numerically for several
$\vv{p}$'s.
Recalling the path integral approach,  the sum of $M$-tuple
$\vv{\ell}$ can be regarded as
a  label of  the gauge equivalent class of flat connections $\alpha$
in~\eqref{expected_asymptotics}, and we can identify the Chern--Simons invariant with
\begin{equation}
  \CS \left( A_{\alpha(\vv{\ell})} \right)
  =
  - \frac{P}{4} \, \left(
    1 + \sum_{j=1}^M \frac{\ell_j}{p_j}
  \right)^2
  \mod 1
\end{equation}
See Refs.~\citen{FintuStern90a,KirkKlas90a} for computations of the
Chern--Simons invariant for the Seifert homology spheres.
Note that this value
originates
from the $T$-matrix~\eqref{T_matrix} of the vector modular form.
Correspondingly, the Reidemeister torsion is given by
\begin{align}
  \sqrt{T_{\alpha(\ell)}}
  & =
  \left|
    \prod_{j=1}^M
    \sin
    \left(
      P \, \frac{\ell_j}{p_j^2} \, \pi
    \right)
  \right| \cdot
  C_{\vv{p}} \left( \vv{\ell} \right)
\end{align}
Here the product of sin-functions originates from the $S$-matrix of
the vector modular form.

\subsection{Ohtsuki Series}

A tail part in the asymptotic formula~\eqref{exact_asymptotics_large_1} has a
simple generating function as was studied in~\eqref{Lawrence_find_Ohtsuki}.

\begin{theorem}
  Let the $T$-series be defined by~\eqref{define_T_large_1}.
  Then the generating function of the $T$-series is
  \begin{equation}
    \label{generate_Tk}
    \frac{
      \prod_{j=1}^M
      \sinh
      \left(\frac{P}{p_j } \, x \right)
    }{
      \left[
        \sinh
        \left( P x        \right)
      \right]^{M-2}
    }
    =
    \frac{1}{2}
    \sum_{k=0}^\infty
    \frac{T_{\vv{p}}(k)}{(2 \, k)!} \, x^{2 k}
  \end{equation}
\end{theorem}
\begin{proof}
  We first study a case of
  $\sum_j \frac{1}{p_j} <1$.
  Using~\eqref{product_and_chi},
  we have
  \begin{equation*}
    \frac{
      \prod_{j=1}^M
      \sinh
      \left(  \frac{P}{p_j } \, x \right)
    }{
      \left[
        \sinh
        \left( P x        \right)
      \right]^{M-2}
    }
    =
    \frac{(-1)^{M-1}}{2}
    \sum_{k=0}^\infty
    \sum_{n=0}^{2 P -1}
    \begin{pmatrix}
      k+ M -3 \\
      k
    \end{pmatrix}
    \,
    \chi_{2P}^{\vv{E}}(n) \,
    \E^{ -
      \left(
        n+ 2 P
        \left( k + \frac{M-3}{2}\right)
      \right) x
    }
  \end{equation*}
  We equate this expression with
  $\sum_{k=0}^\infty \frac{T_k}{(2 \, k)!} \, x^{2k}$.
  Applying the Mellin transformation, we get
  \begin{multline*}
    T_k =
    \frac{(-1)^M}{2} \,
    \frac{
      \left( 2 \, P \right)^{2 k}
    }{
      (M-3)!}
    \sum_{n=0}^{2 P-1}
    \chi_{2 P}^{\vv{E}}(n)
    \sum_{j=0}^{M-3}
    \frac{1}{1 + j + 2 \, k} \,
    B_{1 + j +  2 k}
    \left(
      \frac{n+ P \, (M-3)
      }{
        2 \, P
      }
    \right)
    \\
    \times
    \sum_{k=j}^{M-3}
    S_{M-3}^{(k)} \,
    \begin{pmatrix}
      k \\
      j
    \end{pmatrix} \,
    \left(
      \frac{P \, (M-3) - n}{2 \, P}
    \right)^{k-j}
  \end{multline*}
  Here we have used~\eqref{Hurwitz_Bernoulli}.
  Identities~\eqref{define_poly} and~\eqref{recursion_f_poly} give
  \begin{equation*}
    \sum_{k=j}^{M-3} S_{M-3}^{(k)} \,
    \begin{pmatrix}
      k \\ j
    \end{pmatrix}
    \,
    \left(
      \frac{
        P \, (M-3) - n
      }{
        2 \, P
      }
    \right)^{k-j}
    =
    (-1)^{M+1+j} \,
    (j+1) \,
    f_{j+1}^{M-2}
    \left(
      \frac{n}{2 \, P} - \frac{1}{2}
    \right)
  \end{equation*}
  then we have
  \begin{multline}
    \label{compute_Mellin_Tk}
    T_k
    =
    \frac{1}{2} \,
    \frac{( 2 \, P)^{2k}}{
      (M-3)!
    }
    \sum_{n=0}^{2 P -1}
    \chi_{2 P}^{\vv{E}}(n)
    \\
    \times
    \sum_{j=1}^{M-2} (-1)^j \, \frac{j}{j+2 \, k} \,
    f_j^{M-2}\left(\frac{n}{2 \, P} - \frac{1}{2}\right) \,
    B_{2k+j}\left(
      \frac{n+P\,(M-3)}{2 \, P} 
    \right)
  \end{multline}
  For this expression, we substitute an identity
  \begin{multline*}
    B_{2k+j}\left(\frac{n+P\,(M-3)}{2 \, P}\right)
    =
    B_{2k+j}\left(\frac{n}{2 \, P} + \frac{\eveodd(M)-1}{2}\right)
    \\
    + (2 \, k +j) \,
    \sum_{m=0}^{
      \lfloor \frac{M-4}{2} \rfloor
    }
    \left(
      \frac{n}{2 \, P} + m + \frac{\eveodd(M)-1}{2}
    \right)^{2k+j-1}
  \end{multline*}
  which follows from~\eqref{Bernoulli_difference}.
  As we have
  \begin{multline*}
    \sum_{j=1}^{M-2} j \,
    f_j^{M-2}\left(
      \frac{n}{2 \, P} - \frac{1}{2}
    \right)
    \,
    \left(
      -\frac{n}{2 \, P} - m - \frac{\eveodd(M)-1}{2}
    \right)^{j-1}
    \\
    =
    F_{M-2}\left(
      \frac{n}{2 \, P} - \frac{1}{2},
      - \frac{n}{2 \, P} -m -\frac{\eveodd(M)-1}{2}
    \right)
    = 0
  \end{multline*}
  from~\eqref{generate_f_poly},
  we obtain~\eqref{generate_Tk}.

  When $\sum_j \frac{1}{p_j}>1$, we have another term coming from the
  second term in~\eqref{product_and_chi_others}.
  This gives an additional term to~\eqref{compute_Mellin_Tk};
  \begin{multline*}
    T_k^{\text{add}}
    =
    \frac{(-1)^M}{2} \, \frac{(2 \, P)^{2k}}{(M-4)!} 
    \sum_{a=1}^\infty \sum_{b=0}^{a-1}
    \sideset{}{^\prime}\sum_{\vv{\eta}(a)}
    \sum_{j=1}^{M-3}
    (-1)^j \frac{j}{2 \, k +j} \, f_j^{M-3}\left(
      \frac{1}{2} \sum_i \frac{\eta_i}{p_i} - b - \frac{1}{2}
    \right)
    \\
    \times
    \left(
      B_{2k+j}\left(
        \frac{M-4}{2} + \frac{1}{2} \sum_i \frac{\eta_i}{p_i} -b
      \right)
      +
      (-1)^j \,
      B_{2k+j}\left(
        \frac{M-2}{2} - \frac{1}{2} \sum_i \frac{\eta_i}{p_i} +b
      \right)
    \right)
  \end{multline*}
  Recalling Lemma~\ref{lemma:f_even_odd} and applying the same method
  with above,  we recover~\eqref{define_T_large_1}.
\end{proof}

We give  explicit forms of some $T$-series as follows;
\begin{align}
  T_{\vv{p}}(0) & =0
  \\[2mm]
  T_{\vv{p}}(1) & = 4 \, P
  \\[2mm]
  T_{\vv{p}}(2) & = 8 \, P^3 \,
  \left(
    2 - M +
    \sum_{j=1}^M \frac{1}{p_j^{~2}}
  \right)
  \\[2mm]
  T_{\vv{p}}(3) & = 4\, P^5 \,
  \left(
    5 \, \left(
      \sum_{j=1}^M \frac{1}{p_j^{~2}} + 2 - M
    \right)^2
    -2 \, 
    \left(
      \sum_{j=1}^M \frac{1}{p_j^{~4}} + 2 - M
    \right)
  \right)
\end{align}

Based
on the exact asymptotic expansion~\eqref{exact_asymptotics_large_1} of
the WRT invariant,
we extract the \emph{tail} part and define  the formal $q$-series
$\tau_\infty(\mathcal{M})$
as a quantum invariant of the Seifert homology sphere
$\mathcal{M}=\Sigma(\vv{p})$ 
by identifying $\exp(2 \, \pi \, \I/N)$ with  $q$;
\begin{equation}
  \label{define_tau_infinity}
  q^{\frac{\phi(\vv{p})}{4} - \frac{1}{2}} \,
  (q-1) \cdot
  \tau_\infty(\mathcal{M})
  =
  \sum_{k=0}^\infty
  \frac{T_{\vv{p}}(k)}{k!} \,
  \left(
    \frac{
      \log q
    }{4 \, P}
  \right)^k
\end{equation}
This invariant $\tau_\infty(\mathcal{M})$ of the formal $q$-series
coincides with the  invariant $\tau_N^{\text{int}}(\mathcal{M})$
defined in~\eqref{tau_res_int}.
Namely we have an integral expression for $\tau_\infty(\mathcal{M})$.

The Ohtsuki series~\cite{Ohtsu95b}
$\lambda_n(\mathcal{M})$
is  defined
from the formal $q$-series $\tau_\infty(\mathcal{M})$
by
\begin{equation}
  \tau_\infty(\mathcal{M})
  =
  \sum_{n=0}^\infty
  \lambda_n(\mathcal{M}) \,
  (q - 1)^n 
\end{equation}
Then the Ohtsuki series $\lambda_n(\mathcal{M})$ for
$\mathcal{M}=\Sigma(\vv{p})$  is computed as follows.
\begin{theorem}
  \label{thm:Ohtsuki}
  The Ohtsuki series  $\lambda_n(\mathcal{M})$ for
  the Seifert homology sphere $\mathcal{M}=\Sigma(\vv{p})$ is written
  in terms of $\vv{p}$ and $\phi(\vv{p})$ defined
  in~\eqref{define_phi_p} as
  \begin{align}
    \lambda_n(\mathcal{M})
    & =
    \frac{2}{(n+1)!} \,
    \left.
      \left[
        \prod_{j=0}^n
        \left(
          \frac{P}{4} \frac{\mathrm{d}^2}{\mathrm{d} x^2}+
          \frac{1}{2} - \frac{\phi(\vv{p})}{4} - j
        \right)
      \right] \,
      G(x)
    \right|_{x=0}
  \end{align}
  where the function $G(x)$ is 
  \begin{equation}
    \label{define_G_x}
    G(x)
    =
    \frac{\prod_{j=1}^M \sinh \left( \frac{x}{p_j} \right)}{
      \left[
        \sinh (  x )
      \right]^{M-2}
    }
  \end{equation}
\end{theorem}
\begin{proof}
  From~\eqref{define_tau_infinity} and~\eqref{log_Stirling}, we have
  \begin{equation*}
    (q-1) \, \tau_\infty(\mathcal{M})
    =
    \sum_{m=0}^\infty
    \begin{pmatrix}
      \frac{1}{2} - \frac{\phi(\vv{p})}{4}
      \\
      m
    \end{pmatrix} \,
    (q-1)^m \,
    \sum_{j=0}^\infty \sum_{k=0}^j
    \frac{S_j^{(k)}}{j!} \, \frac{T_{\vv{p}}(k)}{(4 \, P)^k} \, (q-1)^j
  \end{equation*}
  We then have
  \begin{align*}
    \lambda(\mathcal{M})
    & =
     \sum_{j=0}^n
     \frac{1}{(j+1)!} \,
     \begin{pmatrix}
       \frac{1}{2} - \frac{\phi(\vv{p})}{4}
       \\
       n -j
     \end{pmatrix} 
     \sum_{k=1}^{j+1}
     S_{j+1}^{(k)} \,
     \frac{T_{\vv{p}}(k)}{(4 \, P)^k}
     \nonumber
     \\
     & =
     \frac{1}{(n+1)!}
     \sum_{m=0}^n \sum_{k=1}^{n+1-m}
     \begin{pmatrix}
       m+k \\
       m
     \end{pmatrix}
     \, S_{n+1}^{(m+k)} \,
     \left(
       \frac{1}{2} -\frac{\phi(\vv{p})}{4}
     \right)^m \,
     \frac{T_{\vv{p}}(k)}{(4 \, P)^k}
  \end{align*}
  where in the second equality we have  expanded the binomial coefficient
  in terms of the Stirling number of the first kind
  using~\eqref{generate_Stirling}, then we have
  applied~\eqref{Stirling_sum_sum}.
  Theorem~\ref{thm:Ohtsuki} shows that the function
  $G(x)$~\eqref{define_G_x}
  gives
  \begin{equation*}
    \left.
      \left(
        P \, \frac{\mathrm{d}}{\mathrm{d} x}
      \right)^{2 k} \,
      G(x)
    \right|_{x=0} = \frac{1}{2} \, T_{\vv{p}}(k)
  \end{equation*}
  Substituting this expression and
  recalling~\eqref{generate_Stirling},
  we obtain the required formula.
\end{proof}

Explicit forms of the lowest 3 Ohtsuki series $\lambda_n(\mathcal{M})$
for
$\mathcal{M}=\Sigma(\vv{p})$ are
\begin{align}
  \lambda_0(\mathcal{M})
  & =1
  \\[2mm]
  \lambda_1(\mathcal{M})
  & =
  6 \, \lambda_C(\mathcal{M})
  \label{Ohtsuki_and_Casson}
  \\[2mm]
  \lambda_2(\mathcal{M})
  & =
  \frac{3 \, \left( \phi(\vv{p}) \right)^2
    + 12 \, \phi(\vv{p}) -4}{96}
  -\frac{P}{16} \, \left(
    2 - M+\sum_{j=1}^M \frac{1}{p_j^{~2}}
  \right) \, \left( \phi(\vv{p}) +2 \right)
  \nonumber
  \\
  & \qquad
  + \frac{P^2}{96} \,
  \left(
    5 \, \left(
      2 - M+\sum_{j=1}^M \frac{1}{p_j^{~2}}
    \right)^2
    -2 \, \left(
      2 - M+\sum_{j=1}^M \frac{1}{p_j^{~4}}
    \right)
  \right)
\end{align}
where $\lambda_C(\mathcal{M})$ is
the Casson invariant of the Seifert homology sphere
$\mathcal{M}= \Sigma(\vv{p})$~\cite{FukuMatsSaka90,NeumWahl90a}
\begin{equation}
  \label{Casson_sphere}
  \lambda_C(\mathcal{M})
  =
  -\frac{1}{8} + \frac{1}{24 \, P} \left(
    1+ \sum_{k=1}^M \left(
      \frac{P}{p_k}
    \right)^2
    -(M-2) \, P^2
  \right)
  -
  \frac{1}{2} \sum_{k=1}^M
  s\left(\frac{P}{p_k} , p_k\right)
\end{equation}
The relationship between the Casson invariant $\lambda_C(\mathcal{M})$
and the Ohtsuki
series $\lambda_1(\mathcal{M})$ was first proved  in Ref.~\citen{HMuraka93a}.
See Ref.~\citen{CSato97a} for a computation of $\lambda_2(\mathcal{M})$.

\subsection{Lattice Points}

We have seen that
the asymptotics of the WRT invariant is dominated by the
term~\eqref{WRT_dominating},
which shows that
the number of terms in
the  sum of $\vv{\ell}$  are
at most
$D$ defined in~\eqref{define_Dp}.
Though, as was studied in Refs.~\citen{KHikami04b,KHikami04e} for cases
of $M=3$ and $M=4$,
the function
$C_{\vv{p}}(\vv{\ell})$ may vanish
for some $\vv{\ell}$'s.

\begin{theorem}
  \label{thm:lattice_point}
  We fix $M$-tuple $\vv{p}$ with pairwise coprime positive integers
  $p_j$, and let the function
  $C_{\vv{p}}(\vv{\ell})$ be defined by~\eqref{define_Cpl}
  for $\vv{\ell} \in \mathbb{Z}^M$ satisfying
  $1 \leq \ell_j \leq p_j -1 $.
  Due to~\eqref{C_polynomial_invariant}, we have $D$ independent functions.
  We set $\gamma(\vv{p})$ as the number of $M$-tuples $\vv{\ell}$
  satisfying
  \begin{equation*}
    C_{\vv{p}}(\vv{\ell}) \neq 0 ,
  \end{equation*}
  and $L(\vv{p})$ as the integral lattice points
  $\vv{\ell}\in\mathbb{Z}_{>0}^M$
  inside the $M$-dimensional tetrahedron,
  \begin{equation*}
    0< \sum_{j=1}^M
    \frac{\ell_j}{p_j}
    < 1 
  \end{equation*}
  Then we have
  \begin{equation}
    D- \gamma(\vv{p}) \geq L(\vv{p})
  \end{equation}
\end{theorem}

\begin{proof}
  As a generalization of the function $C_{\vv{p}}(\vv{\ell})$ defined
  in~\eqref{define_Cpl},
  we define $C_{\vv{p}}^k(\vv{\ell})$ for $k \geq 0$ by
  \begin{equation}
    C_{\vv{p}}^k(\vv{\ell})
    =
    \sum_{n=1}^{2 P}
    \chi_{2 P}^{\vv{\ell}}(n) \,
    B_k
    \left(
      \frac{n}{2 \, P}
    \right)
  \end{equation}
  in terms of the Bernoulli polynomials.
  We have
  $C_{\vv{p}}(\vv{\ell}) = C_{\vv{p}}^{M-2}(\vv{\ell})$.
  As a generating function
  $Z_{\vv{p}}^{\vv{\ell}}(t)$ of these polynomials, we define
  \begin{equation}
    Z_{\vv{p}}^{\vv{\ell}}(t)
    =
    \sum_{k=0}^\infty
    \frac{
      t^k}{
      k!}
     \, C_{\vv{p}}^k (\vv{\ell}) ,
  \end{equation}
  Using~\eqref{generate_Bernoulli}, we have
  \begin{equation}
    Z_{\vv{p}}^{\vv{\ell}}(t)
    =
    \frac{
      t}{
      \E^t - 1} \,
    \sum_{n=1}^{2 P}
    \chi_{2 P}^{\vv{\ell}}(n) \,
    \E^{\frac{t}{2 P} n}
  \end{equation}
  
  In the case of $0< \sum_j \frac{\ell_j}{p_j} <1$,
  we have
  \begin{equation}
  \label{chi_ell_generate_small}
    \sum_{n=0}^{2 P}
    \chi_{2 P}^{\vv{\ell}}(n) \, z^n
    =
    - z^P \,
    \prod_{j=1}^M \left(
      z^{P \frac{\ell_j}{p_j}} -       z^{-P \frac{\ell_j}{p_j}}
    \right)
  \end{equation}
  which gives
  \begin{equation*}
    Z_{\vv{p}}^{\vv{\ell}}(t)
    =
    -
    \frac{t}{
      \E^{\frac{t}{2}} -      \E^{- \frac{t}{2}}
    } \,
    \prod_{j=1}^M
    \left(
      \E^{\frac{\ell_j}{2 p_j} t} -
      \E^{- \frac{\ell_j}{2 p_j} t}
    \right) .
  \end{equation*}
  This shows that
  \begin{equation*}
    Z_{\vv{p}}^{\vv{\ell}}(t)
    = - 
    \left( \prod_{j=1}^M \frac{\ell_j}{p_j} \right) \,
    t^M +
    O(t^{M+2})    
  \end{equation*}
  and that
  $C_{\vv{p}}^k(\vv{\ell}) =0$
  for $0 \leq k \leq M-1$.
  So we have
  $C_{\vv{p}}(\vv{\ell}) = 0$ when
  $0 < \sum_j \frac{\ell_j}{p_j} <1$.

  The invariance~\eqref{sigma2_chi}
  of the periodic functions $\chi_{2P}^{\vv{\ell}}(n)$
  proves the statement of the theorem.
\end{proof}

In the case of $\sum_j \frac{\ell_j}{p_j} > 1$, the generating function~\eqref{chi_ell_generate_small} should
be replaced with a formula like~\eqref{product_and_chi_others},
and we do not know whether the function $C_{\vv{p}}(\vv{\ell})$
vanishes.
We conjecture that,
when  $\sum_j \frac{\ell_j}{p_j}>1$,
we have
$C_{\vv{p}}(\ell) = 0$ 
iff $\chi_{2P}^{\vv{\ell}}(n)$ coincides with
$\chi_{2P}^{\vv{\ell^\prime}}(n)$ s.t. 
$\sum_j \frac{\ell_j^\prime}{p_j} <1$.

\begin{conj}
  \label{conj:lattice}
  Under the conditions of Theorem~\ref{thm:lattice_point}, we have
  \begin{equation}
    D- \gamma(\vv{p}) = L(\vv{p})
  \end{equation}
\end{conj}

This conjecture was proved for $M \leq 4$ in
Refs.~\citen{KHikami04b,KHikami04e}.
It states that  the  number of the flat connections which contribute
as~\eqref{WRT_dominating} coincides with  the number of integral lattice
points inside the $M$-dimensional tetrahedron.

\subsection{Ehrhart Polynomial}

Explicit form of the number $L(\vv{p})$ of the lattice points inside
the $M$-dimensional tetrahedron was first computed by
Mordell~\cite{Morde51} for cases of $M=3$ and $M=4$;
\begin{itemize}
\item $M=3$;
  \begin{multline}
    L(p_1, p_2, p_3)
    =
    \frac{1}{4} \left( p_1 - 1 \right)\, \left( p_2 - 1 \right)\,
    \left( p_3 - 1 \right)
    +\frac{1}{12  \, P}
    -
    \frac{1}{4}
    \\
    -
    \frac{P}{12} \,
    \left(
      1-
      \frac{1}{p_1^{~2}}
      -   \frac{1}{p_2^{~2}}
      -  \frac{1}{p_3^{~2}}
    \right)
    -
    s(p_1 \, p_2, p_3)
    -s(p_2 \, p_3, p_1) -   s(p_1 \, p_3, p_2)
  \end{multline}

\item $M=4$,
  \begin{multline}
    L(p_1, p_2, p_3, p_4)
    \\
    =
    \frac{1}{8} \prod_{j=1}^4 (p_j - 1)
    + \frac{3}{8} -  \frac{P}{12} +
    \frac{P}{24} \, 
    \sum_{j=1}^4 \frac{1+  \,p_j}{p_j^{~2} }
    +
    \frac{1}{24 \, P} \,
    \left(
      1 - \sum_{j=1}^4 p_j
    \right)
    \\
    -
    \frac{P}{24}
    \sum_{j \neq  k }^4 \frac{1}{p_j^{~2} \, p_k}
    -
    \frac{1}{2}
    \sum_{j=1}^4
    s\left( \frac{P}{p_j} , p_j \right)
    +
    \frac{1}{2}
    \sum_{j \neq k}^4
    s\left(
      \frac{P}{p_j \, p_k}, p_j
    \right)
  \end{multline}
\end{itemize}
Here $s(b,a)$ is the Dedekind sum~\eqref{Dedekind_sum}.
For higher dimension $M$,
the lattice points $L(\vv{p})$ might be written in terms of 
Zagier's higher-dimensional Dedekind sum~\cite{DZagie73a}, but there
seems to exist no applicable expressions.

Although, there is a useful tool to count  the lattice points
(see, \emph{e.g.}, Ref.~\citen{BeckRobi05Book}).
Let $\mathcal{P}$ be the $M$-dimensional open tetrahedron with
integer vertices,
$(p_1, 0, \dots, 0)$, $(0, p_2, 0, \dots, 0), \dots,
(0, \dots, 0 , p_M)$, and
$(0,\dots, 0)$;
\begin{equation}
  \label{tetrahedron_P}
  \mathcal{P}=
  \left\{
    (\ell_1, \dots, \ell_M) \in \mathbb{Z}^M
    ~\Big|~
    \sum_{j=1}^M
    \frac{\ell_j}{p_j} < 1,
    \ell_k > 0
  \right\}
\end{equation}
Let  $E_{\mathcal{P}}(t )$ denote the number of
lattice points in the dilated tetrahedron $t\, \mathcal{P}$.
So we have
\begin{equation*}
  L(\vv{p}) = E_{\mathcal{P}}(t=1) 
\end{equation*}
In the same manner, we suppose that
$E_{\widebar{\mathcal{P}}}(t)$ denotes the number of
lattice points of the closure of $t \, \mathcal{P}$,
\begin{equation*}
  E_{\widebar{\mathcal{P}}}(t)
  =
  \#
  \left\{
    (m_1, \dots, m_M) \in \mathbb{Z}^M
    ~
    |
    ~
    \sum_{j=1}^M \frac{m_j}{p_j} \leq t ,
    m_k \geq 0
  \right\}
\end{equation*}
These functions, $E_{\mathcal{P}}(t)$ and $E_{\widebar{\mathcal{P}}}(t)$,
become polynomials of $t$~\cite{Ehrha67a},
which are  called the Ehrhart polynomial.
Moreover
we  have the Ehrhart--Macdonald reciprocity formula~\cite{IGMacdo71a,Ehrha67a},
\begin{equation}
  E_{\mathcal{P}}(-t )
  =
  (-1)^M \,  E_{\widebar{\mathcal{P}}}(t)
\end{equation}

In general, the number of lattice points
$E_{\widebar{\mathcal{P}}}(t)$ 
becomes  polynomial of $t$ for arbitrary
polytope $\widebar{\mathcal{P}}$~\cite{Ehrha67a}.
We set the coefficients of the Ehrhart polynomial as
\begin{equation}
  E_{\widebar{\mathcal{P}}}(t)
  = c_M(\mathcal{P}) \, t^M + c_{M-1}(\mathcal{P}) \, t^{M-1} +
  \dots
  + c_0(\mathcal{P})
\end{equation}
It is well known that
$c_M(\mathcal{P})$ is the volume of $\mathcal{P}$,
$c_M(\mathcal{P}) = \Vol(\mathcal{P})$,
$c_{M-1}(\mathcal{P})$ is  a   half of the boundary surface area,
$c_{M-1}(\mathcal{P}) = \frac{1}{2} \, \Vol(\partial \mathcal{P})$.
The coefficient
$c_0(\mathcal{P})$ is the Euler characteristic $\chi(\mathcal{P})$,
and $c_0(\mathcal{P}) =1$ when
$\mathcal{P}$ is the convex polytope.

The first nontrivial coefficient of the Ehrhart polynomial for the
$M$-dimensional tetrahedron is thus $c_{M-2}(\mathcal{P})$.
In our case of the $M$-dimensional tetrahedron with
pairwise coprime integers $p_j$, we have~\cite{MBeck00a,BChen02a}
\begin{multline}
  (M-2)! \cdot c_{M-2}(\mathcal{P})
  \\
  =
  \frac{M}{4} +
  \frac{1}{24 \, P} \, \left(
    2 - \sum_{k=1}^M \left(\frac{P}{p_k} \right)^2
    + 3 \, \left(
      \sum_{k=1}^M \frac{P}{p_k}
    \right)^2
  \right)
  -
  \sum_{k=1}^M 
  s \left(
    \frac{P}{p_k} , p_k
  \right)
\end{multline}
Recalling the Casson invariant
$\lambda_C(\mathcal{M})$~\eqref{Casson_sphere} for
the Seifert homology sphere $\mathcal{M}=\Sigma(\vv{p})$, which is
proportional to the first Ohtsuki series~\eqref{Ohtsuki_and_Casson},
we have the following.

\begin{prop}
  The Casson invariant $\lambda_C(\mathcal{M})$ for
  $\mathcal{M}=\Sigma(\vv{p})$ is related to the first nontrivial
  coefficient of the Ehrhart polynomial for the $M$-dimensional
  tetrahedron $\mathcal{P}$~\eqref{tetrahedron_P};
  \begin{equation}
    \lambda_C(\mathcal{M})
    -
    \frac{(M-2)!}{2} \, c_{M-2}(\mathcal{P})
    =
    -\frac{M+1}{8}
    +\frac{M-2}{24} \, P
    -\frac{P}{8} \, \sum_{1\leq j < k \leq M}
    \frac{1}{p_j \, p_k}
  \end{equation}
\end{prop}

It is remarked that the residue formula for $c_{M-2}(\mathcal{P})$
given in Ref.~\citen{MBeck00a} looks like an
expression~\eqref{result_Rozansky} of the WRT
invariant $\tau_N(\mathcal{M})$ for $\mathcal{M}=\Sigma(\vv{p})$.

\section{Some Examples of  Numerical Experiments}
\label{sec:examples}

We give some  numerical experiments on the asymptotic behavior of the
WRT invariants for the Seifert homology spheres.

\subsection{$\Sigma(2,3,5,7 ,11)$}

For $\vv{p}=(2,3,5,7,11)$, we have
$P=2310$, $D=30$, and $\phi(\vv{p})=\frac{34189}{2310}$.
The bases for $30$-dimensional space is given by
$\vv{\ell}=(1,1,\ell_3,\ell_4, \ell_5)$ with
$1 \leq \ell_3 \leq 2$,
$1 \leq \ell_4 \leq 3$,
$1 \leq \ell_5 \leq 5$.
For all these 5-tuples $\vv{\ell}$, we can check that
$C_{\vv{p}}(\vv{\ell}) \neq 0$, which supports
Conjecture~\ref{conj:lattice}
as we have
$\sum_j \frac{1}{p_j} = \frac{2927}{2310}>1$.

In table~\ref{tab:235711} we give numerical results on the Witten invariant
$Z_N(\mathcal{M})$ for $\mathcal{M}=\Sigma(2,3,5,7,11)$, which is
performed with a help of PARI/GP.
We give both  the exact value $Z_N(\mathcal{M})$ and asymptotic value
$Z_N^{(0)}(\mathcal{M})$ for several $N$s.
They vary much with $N$, and comparing these data we see an agreement.

\begin{table}
  \centering
  \begin{equation*}
    \def\arraystretch{1.03}
    \begin{array}{r||r|r}
      N & \text{exact result for $Z_{N}$} &
      \text{asymptotics $Z_{N}^{(0)}$}
      \\
      \hline
      22 &
      -13.346013 + 17.397906 \, \I
      &
      -12.2403 + 16.7013 \, \I
      \\
      23 &
      -0.57682556 - 0.51108147 \, \I
      & 0.020572 + 0.004140\, \I
      \\
      98 &
      0.93263590 - 0.49655457 \, \I
      & 0.323366 + 0.0057023 \, \I
      \\
      99 &
      22.826764 - 367.89360 \, \I
      & 22.8460 - 365.870 \, \I
      \\
      100 & 464.33437 - 287.59556 \, \I
      & 475.688 - 287.973 \, \I
      \\
      998 & 9.2292110 - 9.3324129 \, \I
      & 10.7013 - 1.60581 \, \I
      \\
      999 & -52995.123 - 87204.076 \, \I
      & -53072.7 - 87187.8 \, \I
      \\
      1000 & 694.74344 + 9181.2935 \, \I
      & 683.369 + 9183.49 \, \I
      \\
      2398 & -64.891808 + 46.620794 \, \I
      & -62.4971 + 47.9275 \, \I
      \\
      2399 & 320910.08 + 27551.395 \, \I
      & 321128. + 27510.1 \, \I
      \\
      2400 & 142206.21 - 1871.8080 \, \I
      & 142145. - 1869.06 \, \I
      \\
      2401 & 214250.48 - 80025.187 \, \I
      & 214270. - 79907.4 \, \I
      \\
      \hline
    \end{array}
  \end{equation*}
  \caption{The WRT invariant $Z_N(\mathcal{M})$  for {$\mathcal{M}=\Sigma(2,3,5,7,11)$}.
  Asymptotic formula for $Z_N^{(0)}(\mathcal{M})$ is from~\eqref{WRT_dominating}.
}
  \label{tab:235711}
\end{table}

\subsection{$\Sigma(3,7,8,11,13,17)$}

For $\vv{p}=(3,7,8,11,13,17)$, we have
$\phi(\vv{p})= \frac{338099}{408408}$ and $D=5040$.
This $D=5040$-dimensional vector space is spanned by
$\vv{\ell}=
(1, 1 \leq \ell_2 \leq 3,
 1 \leq \ell_3 \leq 7,
 1 \leq \ell_4 \leq 5,
 1 \leq \ell_5 \leq 6,
 1 \leq \ell_6 \leq 8)$.
Among these $D=5040$ bases, we check that
$C_{\vv{p}}(\vv{\ell}) = 0$ when
$\vv{\ell}=
(1,1,1,1,1,1)$,
$(1,1,1,1,1,2)$,
$(1,1,1,1,1,3)$,
$(1,1,1,1,2,1)$,
$(1,1,1,1,2,2)$,
$(1,1,1,1,3,1)$,
$(1,1,1,2,1,1)$,
$(1,1,1,2,1,2)$,
$(1,1,1,2,2,1)$,
$(1,1,2,1,1,1)$,
$(1,2,1,1,1,1)$,
which supports Conjecture~\ref{conj:lattice}.

Numerical results on the exact value and the asymptotics of the Witten
invariant for $\Sigma(\vv{p})$
are summarized in Table~\ref{tab:378111317}.
We see an agreement.

\begin{table}
  \centering
  \begin{equation*}
    \begin{array}{r||r|r}
      N & \text{exact value $Z_{N}(\mathcal{M})$} &
      \text{asymptotics $Z_{N}^{(0)}(\mathcal{M})$}
      \\
      \hline
      58 & 365.32895 + 679.07006 \, \I
      & 351.149 + 691.982 \, \I
      \\
      59 & 1331.8460 - 433.95047 \, \I
      & 1358.51 - 437.953 \, \I
      \\
      60 & -944.99493 + 765.34451 \, \I
      & -915.949 + 742.606 \, \I
      \\
      61 & 130.91099 + 2814.5744 \, \I
      & 62.8489 + 2763.93 \, \I
      \\
      118 & -0.8206017 + 61.590246 \, \I 
      & 0.782372 + 60.1248 \, \I
      \\
      119 & 8.1857781 + 13.369868 \, \I
      & 0.0195662 + 0.0062675 \, \I
      \\
      120 & 5259.2853 + 4064.4029 \, \I
      & 5232.38 + 4043.94 \, \I
      \\
      121 & 8733.0140 + 5274.8273 \, \I
      & 8659.21 + 5338.15 \, \I
      \\
      238 & -219.36738 - 1.608943 \, \I
      & -216.499 + 1.53462 \, \I
      \\
      239 & -6151.0562 - 5617.75586 \, \I
      & -6220.64 - 5620.95 \, \I
      \\
      240 & -11.492746 + 6.1192358 \, \I
      & 1.67454 + 2.34920 \, \I
      \\
      241 & -26057.019 - 52201.108 \, \I
      & -25950.5 - 52634.8 \, \I
      \\
      242 & 49736.853 - 46390.033 \, \I
      & 49818.0 - 46337.0 \, \I
      \\
      243 & 189895.62 + 265408.04 \, \I
      & 189029. + 265225. \, \I
      \\
      244 & 3782.8814 - 12474.142 \, \I
      & 3814.35 - 12433.5 \, \I
      \\
      998 & 21039.448 + 18091.568 \, \I
      & 21107.1 + 18191.2 \, \I      
      \\
      999 & -12.505553 + 49.861847 \, \I
      &
      -0.0331549 + 0.0338852 \, \I
      \\
      1000 & 78229.306 - 164203.36 \, \I
      &
      78333.1 -   164618. \, \I
      \\
      \hline
    \end{array}
  \end{equation*}
  \caption{The WRT invariant $Z_N(\mathcal{M})$ for $\mathcal{M}=\Sigma(3,7,8,11,13,17)$.
  Asymptotic formula for $Z_N^{(0)}(\mathcal{M})$ is from~\eqref{WRT_dominating}.
}
  \label{tab:378111317}
\end{table}

\section{Conclusion and Discussion}

We have studied the asymptotic expansion of the SU(2)
WRT invariant $\tau_N(\mathcal{M})$ for the
$M$-exceptional fibered Seifert homology spheres
$\mathcal{M}=\Sigma(\vv{p})$  in $N\to\infty$ number theoretically.
We have found that the invariant can be written in terms of a limiting value of
\emph{fractional} derivative, \emph{i.e.} derivative of the Eichler integral, 
of the vector modular forms with weight
$3/2$ and $1/2$.
This supports a result~\cite{KHabi02a} that the WRT invariant
is a \emph{limiting} value  of the holomorphic function in a limit 
that $q$
tends to the $N$-th root of unity.
By use of the nearly modular property of the Eichler integral, we
have obtained  an asymptotic expansion of $\tau_N(\mathcal{M})$ in
the large $N$ limit.

Although an asymptotic behavior of the WRT invariant was previously studied
in, \emph{e.g.}, Ref.~\citen{LawreRozan99a},
the correspondence between modular forms and the quantum invariants enables us to relate
topological invariants such as the Chern--Simons invariants, the
Reidemeister torsion, and the Casson invariant,
with geometries of  modular forms.
For example,
we have found  that the number of the gauge equivalent classes of
flat connections,
which dominate
the WRT invariant in the large-$N$ limit,
is related to the number of integral lattice points inside
the $M$-dimensional tetrahedron.
{}From this view, we have established that
the Casson invariant for the Seifert homology sphere
has a relationship with the 
first non-trivial coefficient of the Ehrhart polynomial.

In our previous papers~\cite{KHikami05a,KHikami05b}
we have shown that
the WRT invariants for the Seifert manifolds with 3-exceptional  fibers 
coincides with a limiting value of the Ramanujan mock theta functions.
Investigated~\cite{KHikami05c} is a
modular  transformation formula of the newly
proposed mock theta functions  based on explicit form
of the WRT invariants.
This intriguing  correspondence seems to originate from
a result that the integral
expression~\eqref{tau_res_int} of the WRT invariant  has a connection with
the Mordell integral~\cite{LJMorde33a}.
Our results presented here will shed a new light on 
geometric and topological aspects of
modular
forms.

Though we have  studied only the SU(2) invariant,
the Witten partition function~\eqref{define_Witten_partition} can be
defined for arbitrary gauge group,
and an explicit form of the invariant  for the Seifert manifold is
given~\cite{HanseTakat04a}.
Extending the   method of Lawrence and Rozansky,
it is shown~\cite{MMarin05a} that
the partition function can be written in
the integral form which   can be interpreted as the matrix model,
and that it  becomes  a sum of local contributions
from the flat connections.
This fact is recently reinterpreted from the viewpoint of the path
integral by use of non-abelian localization~\cite{BeasWitt05a}.
As it is well known that
the Chern--Simons perturbation theory of the SU($N$) Witten
invariant 
as a $1/N$ expansion can be
interpreted from the string theory
(see, \emph{e.g.}, Ref.~\citen{MMarin05b}),
it will be interesting to investigate the quantum invariants/modular
forms correspondence for the WRT invariant associated with SU($N$)
gauge group as a generalization of the present work.

\section*{Acknowledgments}
This work is supported in part  by Grant-in-Aid for Young Scientists
from the Ministry of Education, Culture, Sports, Science and
Technology of Japan.

\appendix
\section{Special Functions and Identities}
\subsection{Dedekind sum}
The Dedekind sum is defined by (see, \emph{e.g.}, Ref.~\citen{RademGross72})
\begin{equation}
  \label{Dedekind_sum}
  s(b,a) =
  \sum_{k=1}^{a-1}
  \Bigl(\Bigl(
  \frac{k}{a}
  \Bigr)\Bigr)
  \,
  \Bigl(\Bigl(
  \frac{k \, b}{a}
  \Bigr)\Bigr)
\end{equation}
where $((x))$ is the sawtooth function 
\begin{equation*}
  (( x ))
  =
  \begin{cases}
    \displaystyle
    x  - \lfloor x \rfloor - \frac{1}{2}
    &
    \text{when $x \not\in \mathbb{Z}$}
    \\[4mm]
    0
    &
    \text{when $x \in \mathbb{Z}$}
  \end{cases}
\end{equation*}
The Dedekind sum can also be written as
\begin{equation}
  s(b,a)
  =
  \frac{1}{4 \, a} \sum_{k=1}^{a-1}
  \cot
  \left(
    \frac{k}{a} \, \pi 
  \right) \,
  \cot
  \left(
    \frac{k \, b}{a} \, \pi
  \right)
\end{equation}

The Dedekind sum is known to satisfy
the reciprocity formula
\begin{equation}
  s(b,a) + s(a,b)
  =
  -\frac{1}{4} +
  \frac{1}{12} \,
  \left(
    \frac{a}{b} + \frac{b}{a}
    +\frac{1}{a \, b}
  \right)
\end{equation}
We note
\begin{gather}
  s(-b, a) = - s (b,a) \\[2mm]
  s(b,a) = s(c,a) \qquad
  \text{for $b \, c = 1 \mod a$}
\end{gather}

\subsection{Gauss sum}
As a discrete analogue of the Gaussian integral,
we have a formula of the Gauss sum as
\begin{equation}
  \sum_{n=0}^{2 N -1} \E^{-\frac{1}{2 N} n^2 \pi \I}
  =
  \sqrt{2 \, N} \,
  \E^{-\frac{1}{4} \pi \I}
\end{equation}
The reciprocity formula of the Gauss sum follows from the Gauss
integral as (see \emph{e.g.} Refs.~\citen{Chandrasek85,LCJeff92a})
\begin{equation}
  \label{Gauss_reciprocity}
  \sum_{n \mod N}
  \E^{\pi \I
    \frac{M}{N} n^2 +
    2 \pi \I k n
  }
  =
  \sqrt{
    \left|
      \frac{N}{M}
    \right|
  }
  \,
  \E^{\frac{\pi \I}{4} \, \sign(N M)} \,
  \sum_{n \mod M}
  \E^{-\pi \I
    \frac{N}{M}
    (n+k)^2
  }
\end{equation}
where $N, M \in \mathbb{Z}$
and
$N \, k \in \mathbb{Z}$,
and
$N \, M$ is even.

\subsection{Bernoulli Polynomial}
The $n$-th Bernoulli polynomial $B_n(x)$ is defined from the generating
function as
\begin{equation}
  \label{generate_Bernoulli}
  \frac{t \, \E^{x t}}{\E^t -1}
  =
  \sum_{k=0}^\infty
  B_k(x) \, \frac{t^k}{k!}
\end{equation}
Some of them are written as follows;
\begin{align*}
  B_0(x)  & = 1 
  \\[2mm]
  B_1(x) & = x- \frac{1}{2}
  \\[2mm]
  B_2(x) & = x^2 - x + \frac{1}{6}
  \\[2mm]
  B_3(x)
  & = x^3 - \frac{3}{2} \, x^2 + \frac{1}{2} \, x
\end{align*}

These polynomials satisfy the following relations;
\begin{gather}
  B_k(1-x) = (-1)^k \, B_k(x) \\[2mm]
  B_k(x+1) - B_k(x) = k \, x^{k-1}
  \label{Bernoulli_difference}
  \\[2mm]
  \frac{\mathrm{d}}{\mathrm{d} x} B_n(x)
  =
  n \, B_{n-1}(x)
\end{gather}
Note that the Bernoulli function  has the Fourier expansion as
\begin{equation}
  B_k(x - \lfloor x \rfloor)
  =
  - k! \,
  \sum_{
    \substack{
      n \in \mathbb{Z} \\
      n \neq 0
    }}
  \frac{\E^{2 \pi \I n x}}{
    \left(
      2 \, \pi \, \I \, n
    \right)^k}
\end{equation}
and that
\begin{equation}
  \label{Bernoulli_x_y}
  B_n(x + y)
  =
  \sum_{k=0}^n
  \begin{pmatrix}
    n \\ k
  \end{pmatrix}
  \, 
  B_k(x) \, y^{n-k}
\end{equation}

\subsection{Stirling Number}

The Stirling number of the first kind $S_n^{(m)}$
denotes the (signed) number of permutations of $n$ elements which
contain $m$ permutation cycles
(see, \emph{e.g.}, Ref.~\citen{Stanley97Book}).
The generating function of $S_n^{(m)}$ is written as
\begin{equation}
  \label{generate_Stirling}
  \prod_{j=0}^{n-1}
  (x- j)
  =
  \sum_{m=0}^n
  S_n^{(m)} \, x^m
\end{equation}
and $S_n^{(m)} \neq 0$ when $n \geq m \geq 0$.
Another form of the generating function is given by
\begin{equation}
  \label{log_Stirling}
  \frac{
    \left( \log q \right)^m
  }{m!}
  =
  \sum_{n=m}^\infty
  S_n^{(m)} \,
  \frac{
    \left( q- 1 \right)^n
  }{
    n!
  }
\end{equation}

Based on  these generating functions, 
we have the recursion relations of $S_n^{(m)}$ as follows;
\begin{gather}
  S_{n+1}^{(m)}
  =
  S_n^{(m-1)}
  -  n \,
  S_n^{(m)}
  \\[2mm]
  \begin{pmatrix}
    m \\ r
  \end{pmatrix} \,
  S_n^{(m)}
  =
  \sum_{k=m-r}^{n-r}
  \begin{pmatrix}
    n \\ k
  \end{pmatrix} \,
  S_{n-k}^{(r)} \, S_k^{(m-r)}
  \label{Stirling_sum_sum}
\end{gather}


\end{document}